\documentclass[12pt]{article}
\usepackage{amsfonts}
\usepackage{mathrsfs}
\usepackage{latexsym}
\usepackage{amssymb}
\usepackage{amssymb,a4}
\usepackage{amsthm,amsmath,anysize}
\def\a{\alpha}
\def\b{\beta}
\def\l{\lambda}
\def\g{\gamma}
\def\G{\Gamma}
\def\0{\bar{0}}
\def\e{\epsilon}
\def\d{\delta}
\def\g{\mathfrak{g}}
\def\dim{\text{dim}}
\def\h{\mathfrak h}
\def\L{\Lambda}
\newtheorem{lemma}{Lemma}[section]
\newtheorem{theorem}[lemma]{Theorem}

\newtheorem{proposition}[lemma]{Proposition}
\newtheorem{definition}[lemma]{Definition}
\newtheorem{corollary}[lemma]{Corollary}

\hoffset \voffset \oddsidemargin=60pt \evensidemargin=45pt
\title{\bf  Representations  for the  restricted Lie color algebras}
\author{Chaowen Zhang \\College of Sciences,\\
 China University of Mining and Technology, \\Xuzhou 221116, China}
\date{ }

\begin{document}
\maketitle


\smallskip\bigskip

\section{Introduction}

  The main goal of the present paper  is to develop the nonrestricted representation
  theory for the restricted Lie color algebras.  The restricted
Lie color algebras was defined in \cite{bm,bp,bp1,f}. Progresses
have been made concerning both the structure theory and the
representation theory of the restricted Lie color algebras in the
literature(see \cite{bm,bp,bp1,sch,f}). But the study on the
representation theory so far has been solely about the  restricted
case. \par
  Assume $\mathbf
F$ is  an algebraically closed field with  char.$\mathbf{F}=p>3$.
Let $\g=\oplus_{\a\in\G}\g_{\a}$ be a restricted Lie color algebra
over $\mathbf F$, where $\Gamma$ is an abelian group with a
bicharacter $(,)$. Unless indicated,  we allow $\G$ to be infinite
throughout the paper, so we consider both the case $\g$ is infinite
dimensional and the case $\g$ is finite dimensional. The aim of this
paper is primarily about the nonrestricted representations of $\g$.
By introducing the $p$-characters for locally finite simple
$\g$-modules,  we work on both restricted and nonrestricted
$\g$-modules. From Sec. 3 on, we are particularly interested in the
modules for a class of restricted Lie color algebras named algebraic
Lie color algebras.
\par
 The paper is arranged as
follows: In Sec. 2 we give the definition of the restricted Lie
color algebras, and define the $p$-character $\chi$ for the simple
locally finite $\g$-modules. Then we study the $\chi$-reduced
enveloping algebras $u_{\chi}(\g)$ for $\g$, and determine the
PBW-type of basis of $u_{\chi}(\g)$. Unlike the Lie algebra case, a
$p$-character $\chi$ is not a linear function unless
$\chi(\g_{\a})=0$ for all $ \a\in \G$ with $p\a\neq 0$.
\par
     In Sec. 3, we  define the algebraic Lie color algebras. We
     introduce
 the FP triple, and  determine the simplicity for the
induced module associated with a FP triple.  \par In Sec. 4, we
  are concerned with the  applications of the main theorems to the algebraic Lie
  color algebra $\g=\text{cgl}(V)$.     We then
  obtain an analogue of the Kac-Weisfeiler theorem. The second application of the
  main theorem is determining the simplicity of the baby
  Verma module $Z^{\chi}(\l)$ for $\g=\text{cgl}(V)$.  When $\g$ is infinite
  dimensional, we employ the $u_{\chi}(\g)-T$
  method given in \cite{f2,j} to work on the modules in the Category
  $\mathcal{O}$.
    \par In the
  appendix, we define the infinite dimensional algebraic group
  $\text{GL}(\{m_i\},\mathbf F)$ and its Lie algebra. Since there are no
   literature accessible to the author on this subject,
   we give a brief introduction for the reader's convenience. The conclusion in this section
  is needed in the definition of algebraic Lie color algebras, and
  independent of the other results in the paper.
\section{Preliminaries}
\subsection{notions and definitions}
Let $\mathbf F^{\times}$ denote the set $\mathbf F\setminus{0}$, and
let $\mathbf N$ denote the set of all nonnegative integers.
\begin{definition}\cite{bm,sch,f} Let $\Gamma$ be an  abelian group. A
bicharacter on $\Gamma$ is a mapping $(,): \G\times \G\rightarrow
\mathbf{F}^{\times}$ such that
$$(\a,\b)(\b,\a)=1,\\(\a,\b+\gamma)=(\a,\b)(\a,\gamma),$$$$(\a+\b,\gamma)=(\a,\gamma)
(\b,\gamma),$$ for all $\a,\b,\gamma\in \G$.
\end{definition}To avoid any confusion, we denote the zero
element in $\G$ by $\bar 0$.\par
  A Lie color algebra $\g$ is a $\G$-graded
$\mathbf F$-vector space $\g=\oplus_{\a\in \G}\g_{\a}$ with a
bilinear map $[,]: \g\times \g\longrightarrow \g$ satisfying
$[\g_{\a},\g_{\b}]\subseteq \g_{\a+\b}$ for every $\a,\b\in\G$ and
$$ (1)\quad [y,x]=-(\b,\a)[x,y],$$$$(2)\quad
(\gamma,\a)[x,[y,z]]+(\a,\b)[y,[z,x]]+(\b,\gamma)[z,[x,y]]=0$$ for
every $x\in\g_{\a}$, $y\in\g_{\b}$, $z\in\g_{\gamma}$ and
$\a,\b,\gamma \in \G$. (1) is called color skew-symmetry and (2) is
called the Jocobi color identity \cite{sch, f}.\par Let $\G$ be an
abelian group with a bicharacter $(,)$. Let
$A=\oplus_{\a\in\G}A_{\a}$ be an associate $\G$-graded $\mathbf
F$-algebra. We define $[x,y]=xy-(\a,\b)yx$, for any $x\in A_{\a}$,
$y\in A_{\b}$, then $A$ becomes a Lie color algebra. We denote this
Lie color algebra usually by $A^-$.\par Let
$\g=\oplus_{\a\in\G}\g_{\a}$ be a Lie color algebra. A color
subalgebra of $\g$ is a $\G$-graded subspace $\g_1$ which is closed
under the Lie bracket operation $[,]$. i.e., $[\g_1,\g_1]\subseteq
\g_1$. A color ideal of $\g$ is a $\G$-graded  subspace $I$ of $\g$
such that $[\g, I]\subseteq I$.  Note that the color symmetry
implies that $[I,\g]\subseteq I$. We can also define the the Lie
color algebra $\g/I$ by letting $[x+I, y+I]=:[x,y]+I$ for $x,y\in\g$
\cite[Sec. 3]{f}. Similarly one defines the solvable and nilpotent
Lie color algebras \cite{bm,f}. The maximal nilpotent ideal of a Lie
color algebra $\g$ is called the nilradical of $\g$. If it is zero,
then we say that $\g$ is reductive. \par A morphism of Lie color
algebras $f: \g\longrightarrow \h$ is a $\mathbf F$-linear mapping
satisfying:\par (1) $f(\g_{\a})\subseteq \h_{\a}$ for any
$\a\in\G$.\par (2) $f([x,y])=[f(x),f(y)]$ for any $x,y\in\g$.\par

\begin{lemma}\cite[Lemma 1.8]{bm} Let $\g$ be a Lie color algebra.
We take $x\in\g_{\a}$, $y\in\g_{\b}$ and $z\in\g_{\gamma}$ with
$\a\in\G^+$ and $\b\in \G^-$. If $2,3\in \mathbf F^{\times}$. Then
$$[x,x]=0, [[y,y],y]=0, [[y,y],z]=2[y,[y,z]].$$
\end{lemma}
 We assume that
char.$\mathbf F=p>3$ throughout the paper. \par  Let $\G$ be an
abelian group with bicharacter $(,)$. Let $V=\oplus_{\a\in
\G}V_{\a}$ be a $\G$-graded vector space. We then define the general
linear Lie color algebra $\text{gl}(V)=\oplus_{\a\in\G}
\text{gl}(V)_{\a}$, where
$$\text{gl}(V)_{\a}=\{f\in \text{gl}(V)| f(V_{\gamma})\subseteq V_{\gamma+\a}, \text{for all}\quad \gamma\in \G\}.$$
Then the associative algebra $\text{gl}(V)$ becomes a Lie color
algebra with Lie multiplication defined by $[f,g]=fg-(\a,\b)gf$, for
all $f\in \text{gl}(V)_{\a}$, $g\in \text{gl}(V)_{\b}$.\par
\begin{definition}Let $A=\oplus_{\a\in \G}A_{\a}$ be a (not
necessarily associative) $\G$-graded algebra and let $d\in
\text{gl}(A)_{\d}$ satisfy $d(xy)=d(x)y+(\d,\a)xd(y)$ for any
$\a\in\G$, $x\in\g_{\a}$, $y\in\g$, then $d$ is called a derivation
of $A$ of degree $\d$.
\end{definition} Let $\g=\oplus_{\a\in\G}\g_{\a}$ be a Lie color algebra.
 By the Jocobi color identity one sees that $adx$ is a derivation for any
homogeneous $x\in \g$.\par Assume $\G$ is a countable set, and let
$\dim V_{\a}=m_{\a}<\infty$,  for each $\a\in \G$. We define a
subalgebra $\text{cgl}(V)$ of $\text{gl}(V)$ as
$$\text{cgl}(V)=\{f\in gl(V)|f(V_{\a})=0 \quad \text{for all but finitely
many }\quad \a\}.$$  Let $\G=\{\a_i\}_{i\geq 1}$ be given the order
that $\a_i<\a_j$ if $i<j$.  Taking an ordered homogeneous basis
$\{v_i\}_{i\geq 1}$ of $V$ such that $\a_i\leq \a_j$, whenever
$v_i\in V_{\a_i}$, $v_j\in V_{\a_j}$. We define $e_{ij}\in
\text{cgl}(V)$ by $e_{ij} (v_k)=\delta_{jk}v_i$. Then
$\{e_{ij}|i,j\in \mathbf{Z}^+\}$ is a basis of $\text{cgl}(V)$, and
we may identify $\text{cgl}(V)$ with the matrix algebra $$\{
(a_{ij})_{i,j\geq 1}|a_{ij}=0 \quad \text{for all but finitely
many}\quad i,j\}.$$  If $\G$ is finite, then
$\text{cgl}(V)=\text{gl}(V)$. In this case we  denote the finite
dimensional general linear Lie color algebra $\text{gl}(V)$ also by
$\text{gl}(m,\G)$, where $m=\dim V=\sum _{\a\in\G}\dim V_{\a}$.
\par Let $\g=\oplus_{a\in\G}\g_{\a}$ be a Lie color algebra.
A $\g$-module $M$ is a $\G$-graded vector space $M=\oplus
_{\gamma\in\G}M_{\gamma}$ such that $\g_{\a}\cdot M_{\b}\subseteq
M_{\a+\b}$ for every $\a,\b\in\G$, and $$[x,y]\cdot m=x\cdot(y\cdot
m)-(\a,\b)y\cdot(x\cdot m)$$ for every $m\in M$, $x\in\g_{\a}$,
$y\in\g_{\b}$.  We say that the $\g$-module $M$ is locally finite
 if there is a finite dimensional $\G$-graded
$x$-invariant subspace $0\neq M_{x}\subseteq M$,  for any
homogeneous $x\in\g$. We  see from Lemma 2.2  that, $\g$ itself is a
locally finite $\g$-module under the adjoint action and the
assumption $p>3$. A $\g$-module $M$ is called simple if it has no
nontrivial $\G$-graded submodules.
\par Let $\g$ be a Lie color algebra.  A representation of $\g$ is a
Lie color algebra homomorphism $\rho: \g\longrightarrow
\text{gl}(V)$ for some $\G$-graded $\mathbf F$-space
$V=\oplus_{\a\in\G}V_{\a}$. Each representation $\rho$ of $\g$
defines a $\g$-module structure on $V$ and vice versa. We say that
the representation $\rho$ is afforded by the $\g$-module $V$.\par

  The bicharacter $(,)$ induces a partition of $\G$ $$\G=\G_+\cup\G_-,$$where
  $\G_{\pm}=\{\gamma\in\G|(\gamma,\gamma)=\pm 1\}$.  We denote $\g^+=\sum_{\gamma
  \in \G^+}\g_{\gamma}$ and $\g^-=\sum_{\gamma\in\G^-}\g_{\gamma}$. $\g^+$ is
  referred as the even subalgebra of $\g$.\par A Lie color
  algebra is called restricted if for every $\gamma\in\G^+$ there
  exists a map $(,)^{[p]}: \g^+_{\gamma}\rightarrow \g^+_{p\gamma}$ such
  that \cite{f} \par (1) $(x+y)^{[p]}=x^{[p]}+y^{[p]}+\sum^{p-1}_{i=1}s_i(x,y)$, where
  $is_i(x,y)$ is the coefficient of $t^{i-1}$ in the polynomial
  $(ad_{\g}(tx+y))^{p-1}(x)\in \g[t]$;\par
  (2) $(rx)^{[p]}=r^px^{[p]}$;\par (3) $ad_{\g}x^{[p]}=(ad_{\g}x)^p$ for every
  $x,y\in\g_{\gamma}$ and every $r\in \mathbf F$.\par
 Let $A=\oplus_{\a\in\G} A_{\a}$ be a $\G$-graded associative algebra. Then $A^-$ becomes
 a restricted Lie color algebra if we define $x^{[p]}=x^p$, the p-th power,
  for any $x\in\g_{\a}$, $\a\in\G^+$.  So that $\text{gl}(V)$ is restricted
  and $\text{cgl}(V)$ is its restricted subalgebra.\par
The Lie color algebra $\text{cgl}(V)$ has a maximal torus
$H=\{e_{ii}|i\geq 1\}\subseteq \g_{\bar 0}$. Let $\dim
V_{\a}=m_{\a}$. Then $\g_{\bar 0}$ is a restricted Lie algebra
$\oplus_{\a_i\in \G}\text{gl}(m_{\a_i})$.
\par Let $\g$ be a Lie color algebra. Consider the two-sided ideal
$I(\g)$ of the tensor algebra $T(\g)$ of $\g$ over $\mathbf F$
generated by $x\otimes y-(\a,\b)y\otimes x-[x,y]$ for $x\in\g_{\a}$,
$y\in\g_{\b}$, $\a,\b\in\G$. Then $U(\g)= T(\g)\slash I(\g)$ is the
universal enveloping algebra of $\g$\cite{f}. \par Similar to the
Lie algebra case, there is an equivalence between the category of
$\g$-modules and the category of unitary $\G$-graded
$U(\g)$-modules\cite[p.114]{f}.\par Let $A=\oplus_{\a\in\G}A_{\a}$
be an associative $\G$-graded algebra, we define the center of $A$
by $Z(A)=\oplus_{\a\in\G}Z(A)_{\a}$, where
$$Z(A)_{\a}=\{a\in A_{\a}|ab=(\a|\b)ba, \text{for all}\quad b\in
A_{\b}, \b\in\G\}.$$  Let $ q\in \mathbf F^{\times}$ and
$n\in\mathbf N$. We denote
$$[n]_q=\frac{1-q^n}{1-q}=\begin{cases} 1+q+\cdots+q^{n-1},
&\text{if}\quad q\neq 1\\n,&\text{if}\quad q=1.\end{cases}$$ Recall
the quantum binomial coefficient
$Q^i_n=:\frac{[n]!_q}{[i]!_q[n-i]!_q}$. Let $a\in Z(A)_{\a}$ and
$b\in A_{\beta}$. Then by induction, one gets the binomial formula
$$(a+b)^n=\sum^n_{i=0}Q^i_n a^ib^{n-i},$$  where $Q^i_n$  is
defined with $q=(\beta,\a)$.\par
 Let
$\g=\oplus_{\a\in\G}\g_{\a}$ be a Lie color algebra. By induction,
one can show that, in the universal enveloping algebra $U(\g)$ of
$\g$,
$$
(adx)^ky=\sum^k_{i=0}(-1)^i(\a|\b)^iC^i_k x^{k-i}yx^i,
$$
for all $x\in\g_{\a}$,  $y\in\g_{\b}$, $\a,\b\in\G$. If $\g$ is
restricted, then we have $$(x^p-x^{[p]})y=(p\a|\b)y(x^p-x^{p]})$$
and hence $x^p-x^{[p]}\in Z(U(\g))_{p\a}$, for all $x\in\g_{\a}$,
 $y\in\g_{\b}$, $\a\in\G^+$, $\b\in\G$.
\par Let $\g$ be a Lie color algebra and
let $U(\g)$ be its universal enveloping algebra. For each $k\geq 0$,
we define a subspace $U_{(k)}$ of $U(\g)$ as follows: If $k=0$, we
let $U_{(0)}=\mathbf{F}$. If $k\geq 1$, we let
$$U_{(k)}=:\langle \Pi^l_{j=1}x_j|x_j\in \g_{\a},\a\in \G;l\geq
0,l\leq k\rangle.$$ Then $\{U_{(k)}\}_{k\in\mathbf N}$ is a
filtration of $U(\g)$.\par  Let $I$ be a totally ordered set.  We
denote
$$N(I)=\{f:I\longrightarrow \mathbf N|f(i)=0 \quad\text{for all but finitely many}
\quad i\in I\}.$$ For each $s\in N(I)$, we denote $|s|=\sum_{i\in
I}s(i)$. Let $\g$ be a Lie color algebra and let $\{z_i\}_{i\in I}$
be a set of elements in $U(\g)$. For each $r\in N(I)$, we use $z^r$
to denote the ordered product $\Pi_{i\in I} z_i^{r(i)}\in U(\g)$.
\begin{lemma} Let $(e_i)_{i\in I}$
be a totally ordered homogeneous basis of
$\g=\oplus_{\a\in\G}\g_{\a}$, and  let $J$ be a subset of $I$.
Assume  $k\in N(I)$.  Suppose there are two sets of homogeneous
elements $(v_i)_{i\in I}$, $(z_i)_{i\in I}$ in $U(\g)$ having the
following properties: \par (1) If $i\in J$, $e^{k(i)}_i=v_i+z_i$,
$v_i\in U_{(k(i)-1)}$, $z_i\in Z(U(\g))$;\par (2) If $i\in
I\setminus J$, there is a $\theta (i)\in J$ such that
$$e_i^{k(i)}=v_i+z_i+t_ie^{k(i)}_{\theta (i)},$$ where $z_i\in
Z(U(\g)), v_i\in U_{(k(i)-1)}$ and
 $t_i\in \mathbf{F}$. Then
$$ B:=\{z^re^s|r,s\in N(I), s(i)<k(i), \text{for all}\quad i\in
I\}$$ is a basis of $U(\g)$,
\end{lemma}
\begin{proof} We prove by induction that $$B_n=\{z^re^s|\sum_{i\in
I}(r(i)k(i)+s(i)\leq n, s(i)<k(i)\}$$ is a basis of $U_{(n)}$ for
any $n\in \mathbf N$. The case $n=0$ is trivial. Assume the claim is
true for the case $n-1$ with $n\geq 1$. By the PBW theorem,
$U_{(n)}$ is spanned by the set
$$P_n=:\{e^{s_1}_{i_1}e^{s_2}_{i_2}\cdots e^{s_k}_{i_k}|i_1<\dots <i_k,
s_1+\dots +s_k\leq n\}.$$ For any  $\Pi_{i}e^{l_i}_i\in P_n$ with
$\sum l_i=n$, assume $l_i=q_ik(i)+r(i)$, $0\leq r(i)<k(i)$, we get
$$\Pi_{i} e^{l_i}_i=\Pi_{j\in J,i\notin J}\cdots
e^{l_j}_j\cdots e_i^{q_ik(i)+r(i)}\cdots$$$$=\Pi_{j\in J, i\notin
J}\cdots e^{l_j}_j\cdots (v_i+z_i+t_ie_{\theta
(i)}^{k(i)})^{q_i}e^{r(i)}_i\cdots$$$$\equiv \Pi_{j\in J, i\notin
J}\cdots e^{l_j}_j\cdots (z_i+t_ie_{\theta
(i)}^{k(i)})^{q_i}e^{r(i)}_i\cdots (\text{mod} U_{(n-1))}$$ (using
the binomial formula given earlier)
$$\equiv \sum_{i\notin J, \eta_i\leq q_i} c(\Pi_{i\notin
J}z^{\eta_i}_i)(\cdots e_j^{l_j+\sum_{j=\theta
(i)}(q_i-\eta_i)k(i)}\cdots e^{r(i)}_i\cdots) (\text{mod}
U_{(n-1)})$$
$$\equiv \sum_{i\notin J,\eta_i\leq q_i}c(\Pi_{j\in J}z^{q_{i,j}}_j)(\Pi_{i\notin J}
z_i^{\eta_i})(\cdots e^{r(j)}_j\cdots e^{r(i)}_i\cdots)(\text{mod}
U_{(n-1)}),$$ where for each $i\notin J$ and $j\in J$,
$$l_j+\sum_{j=\theta(i)}(q_i-\eta_i)k(i)=q_{ij}k(j)+r(j),$$
$0\leq r(j)<k(j)$, and each coefficient $c$ is resulted from the
applications of the binomial formula. It is easy to see that each
term in the summation lies in  $B_n$, so that $U_{(n)}$ is spanned
by $B_n$.
\par To verify the linear independency of $B_n$, we note that
$$z_i^{r(i)}\equiv
\begin{cases} e^{k(i)r(i)}_i, &\text{if}\quad i\in J\\ (e_i^{k(i)}-t_ie_{\theta
(i)}^{k(i)})^{r(i)}, &\text{if}\quad i\notin J\end{cases}
(\text{mod} U_{(k(i)r(i)-1)}).$$ For  each $z^re^s\in B_n$, we have
$$z^re^s\equiv
 \Pi_{i\notin J, j\in J} \cdots e_j^{k(j)r(j)+s(j)}\cdots
(e_i^{k(i)}-t_ie^{k(i)}_{\theta (i)})^{r(i)}e_i^{s(i)}\cdots
(\text{mod}U_{(n-1)})$$$$\equiv \sum_{i\notin J,\eta (i)\leq
r(i)}c\Pi_{j\in J}\cdots e^{k(j)r(j)+s(j)+\sum_{\theta
(i)=j}\eta(i)r(i)}_j\cdots e_i^{k(i)(r(i)-\eta (i))+s(i)}\cdots
(\text{mod}U_{(n-1)}),$$ where each term in the summation is a basis
element in $P_n$ multiplied by some $ c\in \mathbf{F}^{\times}$.
Note that the term  with the minimal power for all $j\in J$ is
$\Pi_{i\in I}e_i^{k(i)r(i)+s(i)}$, so by basic linear algebra we get
the linear independency of $B_n$.
\end{proof}
\begin{definition}Let $A=\oplus_{\a\in\G}A_{\a}$ be an associative
$\G$-graded  $\mathbf F$-algebra. If $ab=(\a,\b)ba$ for any
$a\in\g_{\a}$, $b\in\g_{\b}$ or, equivalently, $A^-$ is abelian,
then we say that $A$ is a color commutative algebra.
\end{definition}
 Let $\g$ be a restricted Lie color algebra. Suppose $\g=\g^+$
and $(e_i)_{i\in I}$ is a totally ordered homogeneous basis of $\g$.
Let $z_i=e_i^p-e_i^{[p]}\in Z(U(\g))$. Then $R=\mathbf F[z_i,i\in
I]$ is a color commutative polynomial ring with indeterminate's
$z_i$.

\subsection{p-characters for the simple modules}
Let $\g=\oplus_{\a\in\G}\g_{\a}$ be a  restricted Lie color algebra.
We say that $\g$ is $p$-finite, if the abelian subalgebra $\langle
x\rangle_p=:\langle x^{[p]},x^{[p]^2},\cdots\rangle$ is finite
dimensional for any homogenous $x\in\g$. Then each finite
dimensional restricted Lie color algebra is $p$-finite. It is also
easy to see that $\g=\text{cgl}(V)$ is $p$-finite. \par Let $\g$ be
a $p$-finite Lie color algebra. In this section we introduce for
each locally finite simple $\g$-module a $p$-character $\chi$.
\begin{proposition}Let $\g=\oplus_{\a\in\G}\g_{\a}$ be a  $p$-finite Lie color
algebra, and $V=\oplus_{\a\in\G}V_{\a}$ be a locally finite simple
$\g$-module. Assume $\a\in\G^+$ such that $p\a$ has finite order
$s$. Then there exists a function $\kappa$:
$\g_{\a}^+\longrightarrow \mathbf F$,  for each $\a\in\G^+$ such
that
$$((x^p-x^{[p]})^s-\kappa(x))\cdot m=0,$$ for any $x\in\g_{\a}$, $m\in
V$.
\end{proposition}
\begin{proof}  Let $x\in\g_{\a}$.  Then
we have $x^p-x^{[p]}\in Z(U(\g))_{p\a}$, and hence
$(x^p-x^{[p]})^s\in Z(U(\g))_{\bar 0}$. Since $V$ is locally finite,
there is a finite dimensional $\G$-graded $x$-invariant subspace
$V_x\subseteq V$. Since $\g$ is $p$-finite, the $\G$-graded subspace
$V'=V_x+\sum_{i\geq 1}x^{[p]^i}V_x$ is finite dimensional and
invariant under the action of both $x$ and $x^{[p]}$.  Acting on
$V'$, the element $(x^p-x^{[p]})^s$ has an eigenvalue $\l$, since
$\mathbf F$ is algebraically closed.

  Let $V_{\l}=\{v\in V|(x^p-x^{[p]})^sv=\l v\}$. It
is easy to see that $V_{\l}$ is a nonzero $\G$-graded submodule of
$V$. The simplicity of $V$ implies that $V_{\l}=V$.  Then the
mapping $x\longrightarrow \l$ defines a function $\kappa (x)$.
\end{proof}
Note that if $\g$ is a finite dimensional restricted Lie color
algebra, then \cite[Lemma 2.5]{bm} says that each simple $\g$-module
is finite dimensional, so that the assumption of the proposition is
satisfied.
\begin{definition}\cite[1.11]{bm} Let $\g=\oplus_{\a\in\G}\g_{\a}$
be a Lie color algebra and let $V=\oplus_{\a\in\G} V_{\a}$ be a
$\g$-module. Suppose $\phi\in \text{gl}(V)_{\a}$.  If $y\phi (
m)=(\b,\a)\phi (y\cdot m)$ for any $y\in\g_{\b}$, $\b\in \G$, $m\in
V$, then we say that $\phi$ is a centralizer of the $\g$-module V.
\end{definition}

 Let $\g=\oplus_{\a\in\G}\g_{\a}$ be a $p$-finite
  Lie color algebra, and let $V=\oplus_{\a\in\G}V_{\a}$ be a simple
$\g$-module as that in the proposition above.  If $\a\in \G^+$ and
$p\a=0$, we let $\chi (x)=\kappa (x)^{\frac{1}{p}}$. Then we have
$$ (x^p-x^{[p]}-\chi(x)^p)\cdot m=0,$$ for any $x\in \g_{\a}$ and
$m\in V$. By \cite[Prop 2.1, p.70]{sf}, the mapping
$x\longrightarrow x^p-x^{[p]}$ from $\g_{\a}$ to $U(\g)$ is
$p$-semilinear. It follows that  $\chi_{|\g_{\a}}\in \g^*_{\a}$. If
$p\a\neq 0$ has  finite order $s$.  Let $x\in \g_{\a}$ with $\kappa
(x)= 0$. Then Proposition 2.5 says that $x^p-x^{[p]}$ is a nilpotent
centralizer of $\g$-module $V$. But $(x^p-x^{[p]})V$ is also a
graded submodule of $V$, so we have $(x^p-x^{[p]})V=0$.
\begin{lemma} If $p\a\neq 0$ and
$\kappa_{|\g_{\a}}\neq 0$, then there exists an invertible
centralizer $\phi\in\text{ gl}(V)_{p\a}$, and $c(x)\in\g^*_{\a}$
such that
$$ (x^p-x^{[p]}-c(x)^p\phi )\cdot m=0$$ for any $x\in\g_{\a}$ and
$m\in V$.\end{lemma}

\begin{proof}  Let $\rho$ be the representation
afforded by the $\g$-module $V$. By assumption, we get $\kappa
(y)\neq 0$ for some $y\in\g_{\a}$.  Then  the centralizer
$\phi=:\rho (y^p-y^{[p]})\in \text{gl}(V)_{p\a}$ is invertible.  For
each $x\in \g_{\a}$, we have
$$\rho(x^p-x^{[p]})\phi^{-1}\in \rho(Z(U(\g))_{\bar 0}).$$ Since the mapping
$x\longrightarrow \rho(x^p-x^{[p]})$ is $p$-semilinear, we get, in
the light of the proof of Prop. 2.6,  a linear function
$c(x)\in\g^*_{\a}$ such that $$[(x^p-x^{[p]})\phi^{-1}-c(x)^p]\cdot
m=0$$ for any $x\in\g_{\a}$, $m\in V$. Thus we have
$[(x^p-x^{[p]})-c(x)^p\phi]\cdot m=0$.

\end{proof}
  Let $\g=\oplus_{\a\in\G}\g_{\a}$ be a
restricted Lie color algebra, and $V=\oplus_{\a\in\G}V_{\a}$ be a
 simple $\g$-module as above. Assume $\rho$ is the
irreducible representation afforded by the $\g$-module $V$.  Suppose
$\a\in \G^+$ has infinite order. Let $x\in\g_{\a}$. Then since
$\rho(x^p-x^{[p]})$ is a centralizer of $\g$-module $V$, we get
either $(x^p-x^{[p]})\cdot V=0$ or $(x^p-x^{[p]})\cdot V=V$. In the
latter case, $\phi=\rho(x^p-x^{[p]})$ is invertible. In the light of
the proof of Prop. 2.6, we get a linear function $c(y)\in\g^*_{\a}$
such that
$$ (y^p-y^{[p]}-c(y)^p\phi )\cdot m=0$$ for any $y\in\g_{\a}$ and
$m\in V$. In the case $V$ is a locally finite $\g$-module, the
freeness of the subgroup $(a)\subseteq \G$ implies that
$x^p-x^{[p]}$ acts nilpotently on $V$, hence annihilates $V$.\par

\begin{proposition} Let $p\a\neq 0$ and $\kappa_{|\g_{\a}}\neq 0$. Assume
$p\a\in \G^+$ has  finite order $s$. There exist $0\neq
\xi\in\g_{\a}$ and $c(x)\in \g^*_{\a}$ such that for any
$x\in\g_{\a}$ and $m\in V$, we have: \par (1)
$[(\xi^p-\xi^{[p]})^s-1]\cdot m=0$.\par (2)
$\rho(x^p-x^{[p]})=c(x)^p\rho(\xi^p-\xi^{[p]})$. \par (3)
$c(x)^{ps}=\kappa(x)$.
\end{proposition}
\begin{proof} By assumption $\kappa (x_0)\neq 0$ for some
$x_0\in\g_{\a}$. Since $\mathbf F$ is algebraically closed, we get
(1) from the semilinearity of the mapping $x\longrightarrow
x^p-x^{[p]}$ from $\g$ into $U(\g)$ and Prop. 2.6.\par (2) We fix
$\xi\in\g_{\a}$ as in (1). Let $\phi=\rho(\xi^p-\xi^{[p]})$. Then
Lemma 2.8 gives  $c(x)\in\g^*_{\a}$ as required.\par (3) By taking
the $s$-th power on both sides of (2) and using Prop. 2.6, one gets
$c(x)^{ps}=\kappa(x)$ for any $x\in \g_{\a}$.
\end{proof}
 Let $\a\in \G^+$  such
that $p\a\neq 0$ has finite order $s$. We define an equivalence
relation on the set
$$\{(\xi,c(x))|c(x)\in \g^*_{\a},\xi\in\g_{\a}, c(\xi)=1\}$$ as
follows: $$(\xi,c(x))\backsim (\eta, b(x))\quad \text{if}\quad
c(x)=c(\eta)b(x),  b(x)=b(\xi) c(x)$$$$\text{and}\quad
c(\eta)^s=b(\xi)^s=1.$$ Each equivalent class is denoted by
$[\xi,c(x)]$.  We denote by $\mathscr{F}_{\a}$ the set all such
equivalency classes. We shall now define the generalized
$p$-character $\chi$.\begin{definition} Let
$\g=\oplus_{\a\in\G}\g_{\a}$ be a restricted Lie color algebra. A
$p$-character of $\g$ is a set $\{\chi_{\a}\}_{\a\in \G^+}$
satisfying the following properties:\par (1) If $p\a=0$, then
$\chi_{\a}\in\g^*_{\a}$.\par (2) If $\a$ has infinite order, then
$\chi_{\a}=0$.\par (3) If $p\a\neq 0$ has  finite order $s$, then
$\chi_{\a}\in \mathscr{F}_{\a}\cup 0$.\end{definition}We denote the
set $\{\chi_{\a}|a\in\G^+\}$ usually  by $\chi$.
\begin{definition}  Let
$\g=\oplus_{\a\in\G}\g_{\a}$ be a restricted Lie color algebra. A
$\g$-module $V=\oplus_{\a\in \G} V_{\a}$ is said to have a
$p$-character $\chi$ if the following statements hold.\par (1) If
$\chi_{\a}\in \g^*_{\a}$, then $(x^p-x^{[p]}-\chi_{\a}(x)^p)\cdot
m=0$ for any $x\in\g_{\a}$, $m\in V$.\par (2) If
$\chi_{\a}=[\xi,c(x)]\in\mathscr{F}_{\a}$, then Prop. 2.9(1),(2) are
satisfied.
\end{definition} It is easy to check that (2) is independent of the
representatives for $\chi_{\a}$.\par  Let $\g$ be a $p$-finite Lie
color algebra. Then  each locally finite simple $\g$-module has a
 $p$-character.
\par Let $\g$ be a restricted Lie color algebra, and let $\chi$ be
a $p$-character for $\g$.  We shall now define the $\chi$-reduced
universal enveloping algebra for $\g$. Let $U(\g)$ be the universal
enveloping algebra of $\g$.   We define a two-sided ideal $I_{\chi}$
of $U(\g)$ generated by the homogeneous elements as follows:$$
x^p-x^{[p]}-\chi (x)^p\cdot 1, \quad x\in \g_{\a},\quad \text{if}
\quad\chi_{\a}\in \g^*_{\a}
$$
$$ (\xi^p-\xi^{[p]})^s-1,
(x^p-x^{[p]})-c(x)^p(\xi^p-\xi^{[p]}),\quad x\in\g_{\a},\quad
\text{if}\quad \chi_{\a}= [\xi, c(x)]\in\mathscr{F}_{\a}.$$  It is
easy to check that the definition of the quotient $U(\g)/I_{\chi}$
is independent of the representatives of each
$\chi_{\a}=[\xi,c(x)]\in\mathscr{F}_{\a}$.\par Let
$A=\oplus_{\a\in\G} A_{\a}$ and $B=\oplus_{\a\in\G} B_{\a}$ be two
$\G$-graded associative algebras. By a homomorphism $f$ from $A$ to
$B$ we usually mean that $f$ is a homomorphism of algebras
satisfying $f(A_{\a})\subseteq B_{\a}$.
\begin{definition}Let $\g=\oplus_{\a\in\G}\g_{\a}$ be a restricted Lie color algebra, and
let $\chi$  be a $p$-character of $\g$. A pair
$(u_{\chi}(\g),\sigma)$ consisting of an associate $\G$-graded
$\mathbf{F}$-algebra with unity and a homomorphism $\sigma:
\g\longrightarrow u_{\chi}(\g)^-$ is called a $\chi$-reduced
universal enveloping algebra if:
\par
(1)$$\quad\sigma (x)^p-\sigma (x^{[p]})=\chi(x)^p, \quad\text{for
any}\quad x\in\g_{\a},\quad \text{if}
 \quad \chi_a\in\g^*_{\a};$$ $$(\sigma(\xi)^p-\sigma(\xi^{[p]}))^s=1, \sigma (x^p)-\sigma
 (x^{[p]})
 =c(x)^p(\sigma(\xi)^p-\sigma(\xi^{[p]}),$$$$ \text{for any} \quad x\in
 \g_{\a},
  \quad\text{if}\quad \chi_{\a}=[\xi,c(x)]\in\mathscr{F}_{\a}. $$
  (2) Given any $\G$-graded associate $\mathbf{F}$-algebra $A=\oplus_{\a\in\G}A_{\a}$ with
  unity and any homomorphism $g:\g\longrightarrow A^-$ such that
  condition (1) is satisfied, there is a unique homomorphism
  $\tilde{g}$: $u_{\chi}(\g)\longrightarrow A$ of associate
  $\G$-graded algebras  such that
  $\tilde{g}\cdot \sigma=g$.
\end{definition}
 Let $\g=\oplus_{\a\in\G}\g_{\a}$ be a restricted Lie
color algebra. Assume $\{e_i\}_{i\in I}$(resp. $\{f_i\}_{i\in M}$)
is a totally ordered homogeneous basis of $\g^+$(resp. $\g^-$). Let
$\chi$ be a $p$-character of $\g$. With a minor adjustment,  we may
choose a basis $\{e_i\}_{i\in I}$ of $\g^+$ containing each $\xi$,
if $\chi_{\a}=[\xi,c(x)]\in \mathscr{F}_{\a}$ for some $\a\in \G^+$.
Let $J$ denote the subset of $I$ $$\{i\in I|e_i=\xi,
\quad\text{if}\quad \chi_{\a}=[\xi,c(x)]\in \mathscr{F}_{\a},
\xi\in\g_{\a}\}.$$ For $\a\in\G^+$,  assume  $p\a\neq 0$ has finite
order $s$.  $s$ varies with different $\a$. We denote all these
different orders simply by $s$.\par Recall the notation $N(I)$ for
each set of indices $I$. We have
\begin{theorem} Let $\g=\oplus_{\a\in\G}\g_{\a}$ be a restricted Lie
color algebra. Assume $\{e_i\}_{i\in I}$(resp. $\{f_i\}_{i\in M}$)
be a totally ordered homogeneous basis of $\g^+$(resp. $\g^-$). Let
$\chi$ be a $p$-character of $\g$. Then the quotient
$U(\g)/I_{\chi}$ is the $\chi$-reduced universal enveloping algebra
having the basis $$\{ \sigma (f)^{\d}\sigma (e)^n|\d\in N(M), n\in
N(I), \d_j=0,1,j\in M; 0\leq n(i)<\begin{cases}p, &\text{if}\quad
i\in I\setminus J\\ ps, &\text{if}\quad i\in J\end{cases} \}.$$
\end{theorem}\begin{proof} The
universal property is obvious.  We now apply Lemma 2.4 to get a
basis of $U(\g)/I_{\chi}$. For each $i\in J$, we choose $k(i)=ps$,
then we have $e_i^{ps}=v_i+z_i$, where
$$v_i=e_i^{ps}-(e^p_i-e_i^{[p]})^s+1\in U_{(ps-1)}$$ and
$$z_i=(e_i^{p}-e_i^{[p]})^s-1\in Z(U(\g)).$$ Both $v_i$ and $z_i$ are homogeneous with
the $\G$-grading $\bar 0$. If $i\notin J$ and $e_i\in \g_{\a}$ with
$\chi_{\a}=[\xi,c(x)]\in\mathscr{F}_a$, we choose $k_i=p$, then we
have $e_i^p=z_i+v_i+t_i\xi^p$, where
$$z_i=e^p_i-e_i^{[p]}-c(e_i)^p(\xi^p-\xi^{[p]})\in Z(U(\g)),$$
$$v_i=e_i^{[p]}-c(e_i)^p\xi^{[p]}\in U_{(p-1)}\quad \text{and}\quad t_i=-c(e_i)^p.$$ If
$i\notin J$ and $e_i\in\g_{\a}$ with $\chi_a\in\g^*_{\a}$, we choose
$k(i)=p$, then we get $e_i^p=v_i+z_i$, where $v_i=e_i^{[p]}\in
U_{(p-1)}$ and $z_i=e_i^p-e_i^{[p]}\in Z(U(\g))$.\par For each $i\in
M$, let  $f_i\in\g_{\a}$,  $\a\in\G^-$. Then since $(\a,\a)=-1$ and
$p>3$, we have in $U(\g)$ that $f^2_i=\frac{1}{2}[f_i,f_i]\in \g$,
so that we may choose $k(i)=2$, $$v_i=\frac{1}{2}[f_i,f_i]\in
\g=U_{(1)},\quad z_i=0.$$ Then we have $f^{k(i)}_i=v_i+z_i$. \par
Applying Lemma 2.4,  we get $$I_{\chi}=\langle z^rf^{\d}e^n|r,n\in
N(I), |r|\geq 1, n(i)<k(i)\rangle,$$ so that $U(\g)/I_{\chi}$ has
the basis claimed, where $\sigma$ is the canonical epimorphism:
$U(\g)\longrightarrow U(\g)/I_{\chi}$.
\end{proof}
We will denote $U(\g)/I_{\chi}$ by $u_{\chi}(\g)$ in the following.
In particular, if $\chi=0$, we simply denote $u_0(\g)$ by $u(\g)$,
which is also referred to as the restricted universal enveloping
algebra\cite[5.7]{f}. Then each locally finite simple $\g$-module is
a simple $u_{\chi}(\g)$-module for some $p$-character $\chi$. A
$u(\g)$-module is also referred to as a restricted $\g$-module.\par

\subsection{The Frobenious algebra $u_{\chi}(\g)$ }
\begin{definition}(1) Let $V=\oplus_{\a\in\G} V_{\a}$ be a $\G$-graded space, and let $f$ be a
bilinear form on $V$. Then $f$ is called symmetric if
$f(x,y)=(\a,\b)f(y,x)$ for any $x\in V_{\a}$, $y\in V_{\b}$.\par (2)
Let $A=\oplus_{\a\in\G} A_{\a}$ be a $\G$-graded algebra. A bilinear
form $f$ on $A$  is called invariant if $f(xy,z)=f(x,yz)$ for any
$x\in A_{\a}, y\in A_{\beta}, z\in A_{\gamma}$, $\a,\b,\gamma\in
\G$.
\end{definition}
\begin{definition}An associate $\G$-graded algebra is called
Frobenious if it has a nondegenerate invariant bilinear form.
\end{definition}
Let $\g$ be a restricted finite dimensional Lie color algebra. Let
$\chi=\{\chi_a\}_{a\in \G^+}$ be a $p$-character of $\g$. We now
show that the  $\chi$-reduced universal enveloping algebra
$u_{\chi}(\g)$ is a Frobenious algebra.\par Taking for $\g^+$(resp.
$\g^-$) a homogeneous basis $\{e_i|i\leq k\}$(resp. $\{e_i|k+1\leq
i\leq n\}$) containing each $\xi$ if
$\chi_{\a}=[\xi,c(x)]\in\mathscr{F}_{\a}$, we denote
$$J=:\{1\leq i\leq k|e_i=\xi \quad\text{for some}\quad
\chi_{\a}=[\xi,c(x)]\}.$$Note that $J=\phi$ simply means that each
$\chi_{\a}$ is a linear function on $\g_{\a}$, such that
$\chi_{\a}=0$ for all $\a\in\G$ with $p\a\neq 0$. \par By the PBW
theorem, $U(\g)$ is a free left $U(\g^+)$-module with a basis
$$\{\Pi^n_{i=k+1}e_i^{a_i}|0\leq a_i\leq 1\}.$$ Let $x\in\g_{\a}$, $\a\in
\G$. For the simplicity, we also  use $\bar x$ to denote $\a$ in the
following. \par Let $R$ be the $\G$-graded subalgebra of $U(\g)$
generated by
$$z_i=\begin{cases} e_i^p-e_i^{[p]}-\chi (e_i)^p\cdot 1,&\text{if}\quad i\leq k,
\quad\text{and}\quad\chi_{\bar e_i}\quad\text{is
linear}\\e_i^p-e_i^{[p]}-c(e_i)^p(\xi^p-\xi^{[p]}),&\text{if}\quad
\chi_{\bar e_i}=[\xi,c(x)], i\notin
J\\(e_i^p-e_i^{[p]})^s-1,&\text{if}\quad i\in J.\end{cases}$$ By
Lemma 2.4, $R$ is a color commutative polynomial algebra with
indeterminate $z_i's$. Then $U(\g^+)$ is free over $R$ with a basis
$$\{\Pi^k_{i=1}e_i^{a_i}|0\leq a_i\leq p-1,\text{if}\quad i\notin J,
0\leq a_i\leq ps-1\quad \text{if}\quad i\in J\}.$$ Recall that we
denote the set of all the mappings from $\{1,\dots,n\}$ to $ \mathbf
N$ by $N(\{1,\dots,n\})$. We now define $\tau\in N(\{1,\dots,n\})$
by
$$
\tau (i)=\begin{cases} p-1,&\text{if}\quad i\leq k \quad\text{and}\quad i\notin J\\
ps-1,&\text{if }\quad i\in J\\ 1,&\text{ if}\quad i\geq
k+1.\end{cases}$$ For two mappings $a,b\in N(\{1,\dots,n\})$, we
define $a\leq b$ if $a(i)\leq b(i)$ for each $1\leq i\leq n$, and
$a<b$ if $a\leq b$ and at least there is $1\leq i\leq n$ such that
$a(i)<b(i)$. Then by Lemma 2.4 $U(\g)$ is free over $R$ with a basis
$\{e^a|0\leq a\leq \tau\}$. As a $R$-module, $$U(\g)=\sum_{a<\tau}
Re^a\oplus Re^{\tau}.$$ We denote the $R$-submodule
$\sum_{a<\tau}Re^a$ by $V$.  We define the $R$-linear map $p_{\tau}:
U(\g)\longrightarrow R$ by $p_{\tau}(v+r e^{\tau})=r$, $v\in V, r\in
R$. Then we get a $R$-bilinear form $\mu$: $U(\g)\times
U(\g)\longrightarrow R$, $\mu(x,y)=p_{\tau}(xy)$, for any
homogeneous $x,y\in \g$. Obviously we have
$\mu(xy,z)=\mu(x,yz)$.\par  Let $\g$ be a Lie color algebra. For
each homogeneous $x\in \g$, the derivation $ad x$ is naturally
extended to a derivation on $U(\g)$ which is also denoted by $ad x$.
A straightforward computation shows that $$ad
x(\Pi^n_{i=1}y_i)=\sum^n_{i=1}(\bar x,\sum_{l=1}^{i-1}\bar
y_l)y_1\cdots y_{i-1}[x,y_i]y_{i+1}\cdots y_n,$$ for any homogeneous
$y_1,\dots,y_n\in\g$.\par
\begin{lemma} Let $\mu$ be the bilinear form defined above,
and let $I$ be a $\G$-graded ideal of $R$. Then $$IU(\g)=\{u\in
U(\g)|\mu(u,U(\g))\subseteq I\}.$$\end{lemma}\begin{proof}(see
\cite[Th.4.2(1)]{sf} for the Lie algebra case) If $J=\phi$, then the
proof for \cite[Th. 4.2]{sf} also works here. Now we assume $J\neq
\phi$. Using the identities given in the proof of Th. 2.13, one can
show that if $a(i)>\tau(i)$ for some $i\in J$, then $e^a\in
RU_{(|a|-1)}$; if $a(j)\leq \tau(j)$ for all $j\in J$, and
$a(i)>\tau(i)$ for some $i\notin J$, then $e^a\equiv c_be^b
(\text{mod} U_{(|a|-1)})$, where $b\leq \tau, |b|=|a|$, $ c_b\in
\mathbf F^{\times}$. In the light of the proof of
\cite[Th.4.2(1)]{sf}, one gets $V=RU_{(|\tau|-1)}$.\par Taking
$a,b\leq \tau$ with $|a+b|\leq \tau$ and using the discussions
above, we get $e^ae^b\equiv
\delta_{a,\tau-b}r(a,b)e^{\tau}(\text{mod} V),$ where $$0\neq
r(a,b)=\Pi^{n-1}_{i=1}(\sum^n_{j=i+1}a_j\bar{e_j}|b_i\bar{e_i}).$$
Then the arguments leading to  \cite[Th 4.2 (1)]{sf} can be repeated
to give the desired result.
\end{proof}
\begin{theorem}Let $\g$ be a restricted Lie color algebra, and
let $\chi$ be a $p$-character of $\g$. Then  $u_{\chi}(\g)$ is
Frobenious.
\end{theorem} \begin{proof} Recall the ideal $I_{\chi}$ of $U(\g)$ defined
earlier.  With a minor adjustment of the proof for
 \cite[Coro 4.3, p. 218]{sf}, we
get the induced bilinear form $\bar {\mu}$:
$$\begin{cases} u_{\chi}(\g)\times u_{\chi}(\g)\longrightarrow
R/I_{\chi}\cong F\\(\bar x, \bar y)\longrightarrow
\mu(xy)+I_{\chi},\end{cases}$$ which is well defined and
nondegenerate.\end{proof}
\par
\begin{definition} Let $\g$ be a restricted Lie color algebra. For
each $x\in\g_{\a}\subseteq \g^+$, if there is $n\in \mathbf Z^+$
such that    $x^{[p]^n}=0$, then we say that  $\g$ is unipotent.
\end{definition} \begin{lemma}Let $\g$ be a finite dimensional
unipotent Lie color algebra. Then $u(\g)$ is a symmetric
algebra.\end{lemma}
\begin{proof} Let $e_1,\dots,e_n$ be a basis of $\g$ consisting of homogeneous
elements.  For a homogeneous $x\in \g$, suppose that $ad
x(e_j)=\sum^n_{i=1}k_{ij}e_i,\quad 1\leq j\leq n.$  A
straightforward calculation shows that
$$ad x (e^a)=\sum^n_{i=1}(\Pi^{i-1}_{j=1}(\bar x,\bar e_j))\frac{1-(\bar x,
\bar e_i)^{a_i}}{1-(\bar x,\bar e_i)}k_{ii}e^a+\sum_{|b|=|a|,b\neq
a} c_be^b(\text{mod} U_{(|a|-1)}).$$

  Since $\g$ is unipotent, $\g$ is a restricted $\g$-module  under the adjoint action. By \cite[Th. 3.2]{f},
$\g$ acts strictly triangulable on $\g$. Then we can find a
homogeneous basis of $\g$, under which the matrix of $ad x$ is
strictly upper triangular for all homogeneous $x\in \g$. Using the
formula above and the identity $V=RU_{(|\tau|-1)}$,  one can easily
get $ [x, e^{\tau}]\in V$  for any homogeneous $x$.
 Then using the fact $[x,e^a]\in U_{(|a|)}$ we obtain  $xu-(\bar x|\bar u)ux\in V$ for all
homogeneous $u\in U(\g)$, and hence $vu-(\bar v|\bar u)uv\in V$ for
all homogeneous $v,u\in U(\g)$. Therefore, $\mu$ is symmetric, so is
the induced bilinear form $\bar\mu$.
\end{proof}
\section{The main theorems}
Let $\g=\oplus_{\a\in\G}\g_{\a}$ be a restricted Lie color algebra.
Let $\text{Aut}(\g)$ be the group of all the automorphisms for $\g$.
We denote
$$\text{Aut}^{res}(\g)=\{f\in \text{Aut}(\g)| f(x^{[p]})=f(x)^{[p]}\quad\text{for
any }\quad x\in\g_{\a},\a\in\G^+\}.$$ Recall the Lie color algebra
$\g=\text{cgl}(V)$ in 2.1. Then $\g_{\bar
0}=\oplus_{\a_i\in\G}\text{gl}(m_{\a_i},\mathbf F)$. If we
 denote the linear algebraic group $\Pi
_{\a_i\in\G}\text{GL}(m_{\a_i},\mathbf F)$ by
$\text{GL}(\{m_{\a_i}\},\mathbf F)$, then we get
$$\text{Lie}(\text{GL}(\{m_{\a_i}\},\mathbf F))=\g_{\bar 0}$$(see the Appendix
for the infinite dimensional case). Acting on $\g$ by conjugation,
$\text{GL}(\{m_{\a_i}\},\mathbf F)$ becomes a subgroup of
$\text{Aut}^{res}(\g)$.

\begin{definition}Let $\g=\oplus _{\a\in\G}\g_{\a}$ be a restricted Lie
color subalgebra of $\text{cgl}(V)$. If there is a linear algebraic
group $G\subseteq \text{Aut}^{res}(\g)\cap
\text{GL}(\{m_{\a_i}\},\mathbf F)$ such that
 $\text{Lie}(G)=\g_{\bar 0}$. Then $\g$ is
referred to as an algebraic Lie color algebra.
\end{definition}
Let $\g$ be an algebraic Lie color algebra.  We fix a maximal torus
$T$ of $G$. Then $\g$ has the root space decomposition relative to
$T$: $\g=H\oplus \sum_{\d\in\Phi}\g_{\d}$. Then we get a triangular
decomposition of $\g$: $\g=N^-\oplus H\oplus N^+$, where
$$N^+=\oplus_{\d\in\Phi^+}\g_{\d}\quad  N^-=\oplus_{\d\in\Phi^-}\g_{\d}
\quad H=\text{Lie}(T).$$ For each $\d\in\Phi^+$, we use
$e_{\d}$(resp. $f_{\d}$, $H_{\d}$) to denote the positive root
vector(resp. negative root vector, $[e_{\d},f_{\d}]$).\par
  In this section we study the simple modules for
 the algebraic Lie color algebras.   We keep the
convention that $\chi\in \g^*_{\bar 0}$ implies that $\chi\in\g^*$
and $\chi(\g_{\a})=0$ for all $\bar 0\neq \a\in\G$. \par Let $\G$ be
an abelian group, and let $V=\oplus_{\a\in\G}V_{\a}$ be a
$\G$-graded $\mathbf F$-vector space. We define the $\mathbf
F$-linear vector space $$V_f^*=\{f\in \text{Hom}_{\mathbf F}(V,
\mathbf F)|f(V_{\a})=0\quad\text{for all but finitely many}\quad
\alpha's\}.$$    Then  $V^*_f=V^*$ in case $V$ is finite
dimensional. \par Let $\g=\oplus_{\a\in\G}\g_{\a}$ be a Lie color
algebra. For each $\a\in\G$, we let
$$\g^*_{f,\a}=\{\psi\in\g^*_f|\psi(\g_{\b})=0\quad\text{for
all}\quad\b\neq \a\}.$$ Then $\g^*_f=\oplus_{\a\in\G}\g^*_{f,\a}$ is
also a $\G$-graded $\mathbf F$-vector space.\par  A linear function
$\chi\in\g^*_{f,\bar 0}$ is called in the standard form if
$\chi=\chi_s+\chi_n$, if $\chi_s(N^{\pm})=0$, $\chi_n(H+N^+)=0$, and
$\chi (H_{\d})\neq 0$ implies that $\chi
(\g_{\d})=\chi(\g_{-\d})=0$. We say that $\chi$ is standard
semisimple(resp. nilpotent), if $\chi=\chi_s$(resp.
$\chi=\chi_n$).\par

\subsection{The algebraic Lie color algebra $\g=\text{cgl}(V)$}
Let $\g=\text{cgl}(V)$ and let $G=\text{GL}(\{m_{\a_i}\},\mathbf
F)$.  If $\G$ is infinite(resp. finite), then the linear algebraic
group $G$  has the maximal torus: $$T=\{(a_{ii})|a_{ii}\in \mathbf
F^{\times},i\in \mathbf Z^+ \} (\text{resp}.
T=\{diag(t_1,\dots,t_n)|t_i\in\mathbf F^{\times}\},
 n=\dim V).$$ We have the root space decomposition
$\g=\oplus_{\d\in\Phi}\g_{\d}$ relative to $T$, where
$$\Phi=\{\e_i-\e_j|i\neq j\},\quad \g_{\e_i-\e_j}=\mathbf Fe_{ij}.$$
We take the set of positive(resp. negative) roots
$\Phi^+=\{\e_i-\e_j|i<j\}$(resp. $\Phi^-=\{\e_i-\e_j|i>j\}$.  Then
$$\Delta=\{\e_i-\e_{i+1}|i=1,2,\cdots,\}$$ is a set of simple
roots.\par
   For each $\d=\e_i-\e_j\in\Phi^+$, we
denote  $$e_{\d}=e_{ij}\quad f_{\d}=e_{ji}\quad
H_{\d}=\begin{cases}e_{ii}-e_{jj},&\text{if}\quad\bar e_{\d}\in\G^+\\
e_{ii}+e_{jj},&\text{if}\quad\bar e_{\d}\in\G^-.\end{cases}$$
 \par  Suppose $e_{ij}\in
\g_{\a}$ and $\a\in\G^+$. Then  the subalgebra
$<e_{ij},e_{ji},e_{ii}-e_{jj}>$ is isomorphic to the Lie algebra
$sl_2$. If $\a\in\G^-$, then the subalgebra
$<e_{ij},e_{ji},e_{ii}+e_{jj}>$ is isomorphic to the Hensenberg
algebra.  $\g$ has the triangular decomposition: $\g=N^++H+N^-$,
where
$$N^{+}=\sum_{i<j}\mathbf Fe_{ij},\quad N^-=\sum_{i>j}\mathbf
Fe_{ij},\quad H=\sum_{i\geq 1}\mathbf Fe_{ii}.$$ Note that if $\bar
e_{ij}=\a\in\G$, then $\bar e_{ji}=-\a$.\par   We now define a
bilinear form on $\g$ by setting $$b(x,y)=(\a,\b)tr(xy),  \quad
x\in\g_{\a}, y\in\g_{\b}.$$ It is easy to see that $b(x,y)=0$ if
$\bar x+\bar y\neq \bar 0\in\G$.
\begin{lemma} (1) $b(,)$ is nondegenerate, that is, for
$x\in\g_{\a}$, $\a\in\G$, $x=0$ if and only if $b(x,y)=0$ for any
homogeneous $y\in\g$.\par (2) $b(,)$ is invariant, that is,
$b([x,y],z)=b(x,[y,z])$, for any $x\in\g_{\a}$, $y\in\g_{\b}$,
$z\in\g_{\gamma}$.
\end{lemma}\begin{proof} (1) is obvious. \par (2) By definition, we
have
$$b([x,y],z)=(\a+\b,\gamma)tr([x,y]z)=(\a+\b,\gamma)tr(xyz)-(\a+\b,\gamma)(\a,\b)
tr(yxz)$$$$=\delta_{\a+\b,-\gamma}(\a+\b,\gamma)tr(xyz)-(\a+\b,\gamma)(\a,\b)
tr(yxz)$$$$=\delta_{-\a,\b+\gamma}(\a,\b+\gamma)tr(xyz)-(\a,\b+\gamma)(\b,\gamma)tr(xzy)$$$$=
b(x,[y,z]),$$ where the second last equality is given by the
identity $$(\a+\b,\gamma)(\a,\b)=(\a,\b+\gamma)(\b,\gamma).$$
\end{proof}
Note that the restriction of $b(,)$ to $\g_{\bar 0}$ is the usual
trace form, which is again nondegenerate and invariant.\par
\begin{lemma}Let $x\in \g_{\a}$, $y\in\g_{\b}$ and let $g\in G=:\text{GL}(\{m_{\a_i}\},\mathbf F)$. Then
$$b(g\cdot x,y)=b(x,g\cdot y).$$ We say that $b(,)$ is $G$-invariant.
\end{lemma}
The proof is a straightforward computation. We leave it to the
interested reader.\par For each homogeneous $x\in\g$, $b(x,-)$
defines a linear function on $\g$. Moreover, $b(x,-)\in \g^*_f$.
$\g^*_f$ is a $G$-module defined by the coadjoint action: $$g\cdot
\chi (x)=:\chi^g (x)=\chi(g^{-1} xg),$$ for any $g\in G$,
$\chi\in\g^*_f$, $x\in \g$. We define a linear mapping $\theta$:
$\g\longrightarrow \g^*_f$ such that $\theta (x)=b(x,-)$. It is easy
to see that $\theta$ is surjective. Then the $G$-invariancy of
$b(,)$ ensures that $\theta$ is an isomorphism of $G$-modules. Since
the restriction of $b(,)$ to $\g_{\bar 0}$ is also $G$-invariant, we
get $\theta (\g_{\bar 0})=\g_{f,\bar 0}^*$.\par For each $x\in
\g_{\bar 0}$, we say that $x$ is semisimple(resp. nilpotent), if
$gxg^{-1}\in H$(resp. $gxg^{-1}\in N^-$) for some $g\in G$. We claim
that each $x\in\g_{\bar 0}$ has a unique Jordan decomposition:
$x=x_s+x_n$, where $x_s$(resp. $x_n$) is semisimple(resp. nilpotent)
and $[x_s,x_n]=0$. This is well known in the case $\g$ is finite
dimensional. Suppose $\g$ is infinite dimensional. Then there is
$n\in \mathbf Z^+$ such that $x\in \oplus_{i\leq
n}\text{gl}(m_{\a_i})$. Therefore $x$ has a Jordan decomposition in
$\oplus_{i\leq n}\text{gl}(m_{\a_i})$. Since the decomposition is
unique in any finite dimensional Lie subalgebra of $\g_{\bar 0}$,
all these decompositions  agree.  So we can define it to be the
Jordan decomposition for $x\in \g_{\bar 0}$.\par Let
$\chi\in\g^*_{f,\bar 0}$. By our convention this means
$\chi(\g_{\a})=0$ for all $0\neq \a\in\G$. Then there is
$N\in\mathbf Z^+$ such that $\chi (\text{gl}(m_{\a_i}))=0$ for all
$i>N$.  Applying a similar argument as that for finite dimensional
algebraic Lie algebras\cite{k3}, and using the Jordan decomposition
in $\g$, we get $\chi=\chi_s+\chi_n$, and there is $g\in
\oplus^N_{i=1} \text{GL}(m_{\a_i})\subseteq G$, such that
$\chi_s^g(N^{\pm})=0$ and $\chi_n^g(H+N^+)=0$. Besides,
$\chi_s^g(H_{\a})\neq 0$ implies $\chi^g(e_{\a})=\chi^g(f_{\a})=0$,
so that $\chi^g=\chi^g_s+\chi^g_n$ is in the standard form.
\par
Let $\chi\in \g^*_{f,\bar 0}$. For each $g\in G$, obviously the
$\G$-graded algebras  $u_{\chi}(\g)$ and $u_{\chi^g}(\g)$ are Morita
equivalent, so we can assume  that $\chi$ is in the standard form in
the rest of the paper.\par Let $\chi=\chi_s+\chi_n\in\g^*_{\bar
0,f}$ be in the standard form. For each $\a\in\G$, we define
$$\mathcal{Z}_{\a}=:c_{\g_{\a}}(\chi_s)=\{x\in\g_{\a}| \chi_s([x,y])=0\quad
\text{for any homogeneous}\quad y\in\g\}.$$ Then
$\mathcal{Z}=\oplus_{\a\in\G}\mathcal{Z}_{\a}$ is a Lie color
subalgebra of $\g$. In particular,  we have
$$\mathcal{Z}=\oplus_{\d\in\Phi,\chi(H_{\d})=0}\g_{\d}\oplus H.$$ Since
$\chi\in\g^*_{\bar 0,f}$, the codimension of $\mathcal{Z}$ in $\g$
is finite. Let $P= \mathcal{Z}+N^+$. Then $P$ is a parabolic
subalgebra of $\g$. Let $\mathscr N^+$ denote $P^{\perp}$ with
respect to the invariant form $b(,)$. It is easy to check that $$
\mathscr N^+=\oplus_{\d\in\Phi^+,\chi(H_{\d})\neq 0}\g_{\d},$$   and
$ \mathscr N^+$ is  the finite dimensional  nilradical of $P$. \par
\subsection{Category $\mathcal {O}$}

 Let $\g$ be an algebraic Lie color algebra, and let   $\chi=\chi_s+\chi_n\in \g^*_{f,\bar 0}$
 be in the standard form. Recall $G\subseteq \text{Aut}^{res}(\g)$ with a maximal torus $T$ such that
 Lie$(G)=\g_{\bar 0}$.  We choose $T_0$ to be a subgroup of
 $T$ such that $\chi_s(\text{Lie}(T_0))=0$, and $\chi(Ad(t)x)=\chi(x)$
 for all $x\in\g$,  $t\in T_0$.\par
 Let $\mathscr {S}$ denote a restricted Lie color subalgebra
 of $\g$ containing $\text{Lie}(T_0)$.  We identify $u_{\chi}(\mathscr{S})$ with its
canonical image in $u_{\chi}(\g)$. Following \cite{f2,jj}, we define
the
 $u_{\chi}(\mathscr{S})-T_0$ modules. A $\G$-graded $\mathbf
 F$-vector space $V=\oplus_{\a\in\G}V_{\a}$ is called a  $u_{\chi}(\mathscr{S})-T_0$ module
 if it is both a $ u_{\chi}(\mathscr{S})$-module and a $T_0$-module
 such that the following conditions holds:\par (1) We have
 $t(xv)=Ad(t)(x)(tv)$ for all $t\in T_0$, $x\in\mathscr{S}_{\a}$,
 $v\in V_{\b}$, $\a,\b\in \G$.\par (2) The restriction of the
 $\mathscr{S}$-action on $V$ to Lie$(T_0)$ is equal to the
 derivative of the $T_0$-action on $V$.\par

  Let $u_{\chi}(\g)$ be the
$\chi$-reduced universal algebra of $\g$, and let
$u_{\chi}(\g)-\text{Mod}_{T_0}$ be the category of
$u_{\chi}(\g)-T_0$-modules. We define the category $\mathcal {O}$ to
be the full subcategory of $u_{\chi}(\g)-\text{Mod}_{T_0}$ whose
objects satisfying the following conditions(\cite{hh}).\par (1)
$M=\oplus_{\a\in\G}M_{\a}$ is a finitely generated
$u_{\chi}(\g)-T_0$-module.\par (2) $M=\oplus_{\a\in\G}M_{\a}$ is
$H-T_0$-semisimple, that is, $M$ is a weight module:
$M_{\a}=\oplus_{\l\in H^*}M_{\l}=\oplus_{\l\in X(T_0)}M_{\l}$ for
each $\a\in \G$.\par (3) $M$ is locally $N^+$-finite: for each
homogeneous $v\in M$, the $\G$-graded subspace $u(N^+)\cdot v$ of
$M$ is finite dimensional.\par  It is also easy to see that
$\mathcal{O}$
  is closed under submodules, quotients, and finite direct sums.\par
  Example:  Let $\g=\text{cgl}(V)$ and let $\dim
\g=\infty$. Assume $\chi\in\g^*_{f,\bar 0}$.  Then there is $N\in
\mathbf Z^+$ such that $\chi(\text{gl}( m_{\a_i}))=0$ for all $i\geq
N$. Let
 $T_{\a_i}$ denote the maximal torus of $\text{GL}(m_{\a_i},\mathbf F)$ consisting of
 diagonal matrices.
 We choose $$T_0=:\Pi _{i\geq N}T_{\a_i}=\{(A_i)|A_i=E_{\a_i}\quad\text{
 for all}\quad i< N, A_i\in T_{\a_i}\quad\text{for all}\quad i\geq N\}.$$
  Then we see that $\chi(Ad(t)x)=\chi(x)$ for all $
 t\in T_0$ and $x\in \g_{\bar 0}$.\par
  By assumption $\chi (N^+)=0$, $u(N^+)\subseteq u_{\chi}(\g)$ is
  the reduced enveloping algebra for $N^+$. Then it is easy to see that
  $\g=\text{cgl}(V)$, as a $\g$-module under the adjoint action, is in
  Category $\mathcal{O}$ for $\chi=0$ and $T_0=T$.

 \subsection{The Harish-Chandra
homomorphism} Let $\g=\sum_{\d\in\Phi}\g_{\d}$ be an algebraic Lie
color algebra, where $\Phi=\Phi^+\cup\Phi^-$ is its root system.  We
define a mapping $\bar{p}: \cup \g_{\a}-\{ 0\}\longrightarrow
\mathbf N$ by $$\bar p(x)=\begin{cases} p,&\text{if}\quad
x\in\g_{\a}, \a\in\G^+\\2,&\text{if} \quad
x\in\g_{\a},\a\in\G^-.\end{cases}$$ For the simplicity,
 we write $\bar p$ instead of $\bar p(x)$.\par  Let $\chi\in
\g_{f,\bar 0}^*$ be standard semisimple. For each $h\in H$, the
derivation $ad h$ of $\g$ can be extended naturally to the
$\chi$-reduced universal enveloping algebra $u_{\chi}(\g)$. Then
$u_{\chi}(\g)$ has a weight space decomposition:
$u_{\chi}(\g)=\oplus_{\l\in H^*} u_{\chi}(\g)_{\l}$. It is easy to
see that $u_{\chi}(\g)_0$ is spanned by
$$e_{\d_1}^{k_1}\cdots e_{\d_n}^{k_n}u_{\chi}(H)f_{\d_1}^{k_1}\cdots
f_{\d_n}^{k_n}, \d_i\in \Phi^+, 0\leq k_1,\dots, k_n\leq
\bar{p}-1.$$ Recall $N^+=\sum_{\d\in \Phi^+}\g_{\d}$ and
$N^-=\sum_{\d\in\Phi^+}\g_{-\d}$. Applying a similar argument as
that for \cite[7.4]{d}, we have
\begin{lemma}Let $L=u_{\chi}(\g)N_+\cap u_{\chi}(\g)_0$. Then\par
(1) $L=N^-u_{\chi}(\g)\cap u_{\chi}(\g)_0$, and $L$ is a two-sided
ideal of $u_{\chi}(\g)_0$.\par (2) $u_{\chi}(\g)_0=u_{\chi}(H)\oplus
L$.
\end{lemma}
The projection of $u_{\chi}(\g)_0$ onto $u_{\chi}(H)$ with kernel
$L$ is called the $Harish-Chandra$ homomorphism, denoted by
$\gamma$.\subsection{FP triples} \begin{definition}Let $\Phi$ be the
root system of an algebraic Lie color algebra $\g$. We call a subset
$\Phi_1\subseteq \Phi$ a subsystem(resp. an additive subset) of
$\Phi$ if: for any $\d_1,\d_2\in\Phi_1$,
$l_1\d_1+\l_2\d_2\in\Phi_1$(resp. $\d_1+\d_2\in\Phi_1$), whenever
$l_1\d_1+l_2\d_2\in\Phi$, $l_1,l_2\in \mathbf Z$(resp.
$\d_1+\d_2\in\Phi$).\end{definition}
  Let $\g$ be the Lie algebra of a
semisimple, simply connected algebraic group $G$. Let $\chi\in \g^*$
be in the standard form $\chi=\chi_s+\chi_n$. By \cite{f1}, $\chi$
induces a parabolic subalgebra $P=\mathcal Z\oplus N$ of $\g$, where
$\mathcal Z$ is a Levi factor of $P$ and $N$ is the nilradical. Let
$T$ be the maximal torus of of $G$,
 and let $\Phi^+$,
$\Phi_1$ and $\Phi_u$ denote the positive roots of $T$ in $\g$ of
$\mathcal Z$ and $N$. Then Friedlander and Parshall proved that
\cite[Lemma 8.4]{f1}, there is an ordering of $\Phi_u$ as
$\d_1,\dots,\d_t$ such that for each $1\leq i\leq m+1$,
$$\Phi^+_i=\Phi_1^+\cup \{-\d_1,\dots,-\d_{i-1},\d_i,\dots,\d_t\}$$
is a system of positive roots for $\Phi$ in which $\d_i$ is a simple
root. For each $i$, $\{-\d_1,\dots,-\d_i\}$ is an additive subset of
$\Phi$ normalized by $\Phi^+_1$. \par
\begin{definition} Let $\g=\oplus _{\a\in\G}\g_{\a}$ be an algebraic Lie
color algebra with root system $\Phi$, and let $\Phi_1$ be a
subsystem of $\Phi$ such that the set $\Phi\setminus \Phi_1$ is
finite. Put
 $\mathcal Z=\oplus_{\d\in \Phi_1} \mathcal \g_{\d}\oplus H$. Let
$\Phi^+-\Phi^+_1=\{\d_1,\dots,\d_m\}$. For each $i$, suppose
$$\Phi^+_i=: \Phi^+_1\cup \{-\d_1,\dots,-\d_{i-1},\d_i,\dots\d_m\}$$ is a positive
root system of $\g$, in which $\d_i$ is a simple root,  and
 $\{-\d_1,\dots,-\d_i\}$ is an additive subset of $\Phi$
normalized by $\Phi^+_1$. Then we call
 $(\g,\mathcal Z,\{\d_i\}^m_{i\geq 1})$  a
Friedlander-Parshall triple, or simply a FP triple.
\end{definition}

\par Examples: (1) Let $\g=\text{gl}(m,\G)$  and let $\Phi=\Phi^+
\cup\Phi^-$ be its root system. We denote the root system of
$\g_{\bar 0}$ by $\Phi_0$. Let $N_0^+=\sum_{\d\in \Phi_0^+}\g_{\d}$.
 Assume $\chi\in\g^*_{\bar 0}$ is in the standard form:
$\chi=\chi_s+\chi_n$, where $\chi_s$(resp. $\chi_n$) is standard
semisimple(resp. nilpotent). Let $\mathcal Z=c_{\g}(\chi_s)$. Then
$$\mathcal Z=\sum_{\chi_s(H_{\a})=0}\g_{\a}+H, \quad \mathscr
N^+=\sum_{\a\in\Phi^+, \chi_s(H_{\a})\neq 0}\g_{\a}. $$ By
\cite[8.4]{f1} and Lemma 4.1, there is an order of the  roots of
$\mathscr N^+$: $\d_1,\dots,\d_m$, such that $(\g,\mathcal
Z,\{\d_1,\dots,\d_m\})$ is a  FP triple. \par (2) Let
$\g=\text{gl}(m,\G)$  and let $H$ be its maximal torus as above.
 By the argument in the proof of
\cite[8.4]{f1} and Lemma 4.1, we can put $\Phi^+$ in such an order
that $(\g, H,\Phi^+)$ is a FP triple.
\par
(3) Let $\g=\text{cgl}(V)$ be infinite dimensional, and let
$\chi\in\g^*_{\bar 0,f}$ be in the standard form. Recall the  Lie
color subalgebra $\mathcal Z$ and $\mathscr N^+$. We claim that an
order can be given to the roots of $\mathscr N^+$: $\d_1,\dots,
\d_m$, such that $(\g,\mathcal Z,(\d_i)^m_{i=1})$ is a FP triple.
For each  $n\in\mathbf Z^+$, let us denote
$\g_n=:\text{gl}(\sum_{i\leq n} m_{\a_i},\mathbf F)\subseteq \g$.
Then $\g_n$ is a Lie color subalgebra of $\g$ with the $\bar 0$
component $\oplus ^n_{i=1}\text{gl}(m_{\a_i},\mathbf F)$. Moreover,
$(\g_n)^{\infty}_{n=1}$ is a filtration for $\g$:
$$\g_1\subseteq \g_2\subseteq \cdots \subseteq \g=\cup ^{\infty}_{n=1}\g_n.$$
 Let $n$ be large enough
that the root system of $\g_n$ contains $\d_1,\dots,\d_m$. Then the
claim follows from  \cite[Lemma 8.4]{f1} and Lemma 4.1.\par

Let $(\g, \mathcal Z,\{\d_1,\dots,\d_m\}) $ be a FP triple. For each
$0\leq k\leq m$, we denote the subspace of $\g$
$$\mathscr N_k^+=:\sum^{k}_{i=1}\g_{-\d_i}+\sum^m_{i=k+1}
\g_{\d_i}.$$  Let $P_0$ be the subspace of $\g$ defined by
 $$P_0=: \mathcal
Z+\mathscr N^+_0=\mathcal Z+\sum^m_{i=1}\g_{\d_i}.$$   Then we have
\begin{lemma}$P_0$ is a parabolic Lie color subalgebra of $\g$ with $\mathscr N^+_0$
as its  nilpotent radical.
\end{lemma}\begin{proof}Since $\{\d_1,\d_2,\dots,\d_m\}$
is an additive subset of $\Phi^+$, $\mathscr{N}^+_0$ is a nilpotent
subalgebra of $\g$. Let $\d\in\Phi^+_1$. By definition, we have
$$\{\d-\d_1,\d-\d_2,\dots,\d-\d_m\}=\{-\d_1,-\d_2,\dots,-\d_m\},$$ so
that $[f_{\d},\mathscr{N_0^+}]\subseteq \mathscr{N}^+_0$. Let
$\d\in\Phi^+_1$. Suppose $\d+\d_i\notin \{\d_1,\d_2,\dots,\d_m\}$
for some $1\leq i\leq m$. Then we get $\d+\d_i\in \Phi^+_1$, and
hence $\d_i=(\d+\d_i)-\d\in\Phi_1$, contrary to the assumption that
$\Phi_1$ is a subsystem of $\Phi$. Then we conclude that
$[e_{\d},\mathscr{N}^+_0]\subseteq \mathscr{N}^+_0$. It follows that
$[\mathcal Z, \mathscr N^+_0]\subseteq \mathscr N^+_0$. Thus, $P_0$
is a parabolic subalgebra of $\g$ and $\mathscr{N}^+_0$ is an ideal
of $P_0$.
\end{proof}Let us denote
 $P_0^-=: \mathcal
Z+\mathscr N^+_m=\mathcal Z+\sum^m_{i=1}\g_{-\d_i}.$  Similarly one
can prove
 that $P_0^-$ is a parabolic Lie color subalgebra of $\g$
with $\mathscr N^+_m$ as its nilpotent radical.
\subsection{Induced modules associated with the FP triples}
Let $\g$ be an algebraic Lie color algebra, and let
$(\g,\mathcal{Z},(\d_i)^m_{i=1})$ be a FP triple. Assume
$\chi=\chi_s+\chi_n\in\g^*_{f,\bar 0}$ is in the standard form, and
assume
 $\chi(\g_{\pm\d_i})=0$ for all $1\leq i\leq m$. Let $\mathscr M$ be
 a simple $u_{\chi}(P_0)$(resp. $u_{\chi}(P_0^-)$) module. Then by
 \cite[Coro. 3.2]{f}, $\mathscr{M}$ is annihilated by $\mathscr
 N^+_0$(resp. $\mathscr N^+_m$). \par \begin{definition}
  Let $\g=N^-+H+N^+$ be the
triangular decomposition of $\g$, and let $M$ be a $\g$-module.
Assume $0\neq v\in M_{\a}$. If there is $\l\in H^*$ such that
$N^+v=0$ and $hv=\l(h)v$ for all $h\in H$, then  $v$ is called a
maximal vector with $H$-weight $\l$.\end{definition} Let $B=H+N^+$.
Then each maximal vector $v$ defines a 1-dimensional
$u_{\chi}(B)$-module with $\l$ satisfying $\l(h)^p-\l(h)=\chi(h)^p$,
$h\in H$, and $v$ also defines an induced $u_{\chi}(P_0)$-module
$u_{\chi}(P_0)\otimes_{u_{\chi}(B)}\mathbf Fv$. We denote each
simple quotient of $u_{\chi}(P_0)\otimes_{u_{\chi}(B)}\mathbf Fv$ by
$\mathcal M(\l)$.\par Remark: $\mathcal M(\l)$ is a simple
$u_{\chi}(P_0)$-module generated by a maximal vector $v$ of weight
$\l$.
 $v$ may not be unique, so that $\l$ is not necessarily unique. But
 this does no harm to our discussions in the following.  If
 $\chi$ is standard semisimple, one can easily show that
 $u_{\chi}(P_0)\otimes_{u_{\chi}(B)}\mathbf Fv$ has a unique simple
 quotient $\mathcal M(\l)$, and $\mathcal M(\l)$ contains a unique maximal
 vector $v$(with $H$-weight $\l$). \par

 Suppose $\mathcal M(\l)$ is a simple
$u_{\chi}(P_0)$-module as that given above. We now define the
induced module
$$ Z^{\chi}(\mathcal M(\l))=u_{\chi}(\g)\otimes_{u_{\chi}(P_0)}\mathcal M(\l).$$
\par
First, let $\g$ be finite dimensional. Then each simple
$u_{\chi}(\g)$-module $M$ is  finite dimensional.  We have by using
\cite[Th.3.2]{f} that each simple $u_{\chi}(P_0)$-submodule of $M$
contains a maximal vector, and hence in the form  $\mathcal M(\l)$,
so that $M$ is a homomorphic image of $Z^{\chi}(\mathcal
M(\l))$.\par Secondly, let $\g=\text{cgl}(V)$ be infinite
dimensional, and let $(\g,\mathcal{Z},(\d_i)^m_{i=1})$ be the FP
triple  given in Example (3) of Sec 3.4, with $T_0$ given in Example
3.2.  We  now have
\begin{lemma} Each simple $u_{\chi}(P_0)$-module $\mathcal M$ in  Category
$\mathcal{O}$ contains a maximal vector.\end{lemma} \begin{proof}
Recall the definition of $T_0$. Let us denote
$\underline{m}=:\sum_{i<N}m_{\a_i}$. From 3.2, we have
$$T_0=\{(t_{ii})| t_{ii}\in\mathbf F^{\times}, t_{ii}=1\quad\text{for all}
\quad i\leq \underline{m}\}.$$ Let us denote the weight lattice of
$T_0$ by $\L=\oplus^{\infty}_{i=\underline{m}+1}\mathbf Z\e_i.$ For
$\l=\sum n_i\e_i\in \L$, we define the $T_0$-height of $\l$ to be
$|\l|=\sum n_i$. It should be noted that the weight of the positive
root vector $e_{ij}$ relative to $T_0$ is
$$\text{wt}(e_{ij})=\begin{cases} \e_i-\e_j,&\text{if}\quad
\underline{m}<i<j\\-\e_j,&\text{if}\quad i\leq
\underline{m}<j\\0,&\text{if}\quad i<j\leq
\underline{m}.\end{cases}$$

  Let $\mathscr M$
be a simple $u_{\chi}(P_0)$-module in the category $\mathcal{O}$ of
$u_{\chi}(P_0)-T_0$-modules. Assume $0\neq u\in \mathscr M$ is a
homogeneous weight vector. Since $\mathscr M$ is locally
$N^+$-finite, $\mathcal N=:u(N^+) u\subseteq \mathscr M$ is a finite
dimensional $u_{\chi}(B_0)-T_0$-submodule.   Recall the $T_0$-action
on $\mathcal N$ is given by $$ t(xu)=Adt(x)tu.$$Then using the PBW
theorem, with a given order for the positive root vector $e_{ij}'s$,
$\mathcal N$ is spanned by a finite set of $H-T_0$-weight vectors:
$$\mathcal S=:\{\Pi_{i<j} e_{ij}^{k_{ij}} u|0\leq k_{ij}\leq \bar p\}.$$ Assume
the $T_0$-weight of $u$ is $\l\in X(T_0)$. Then each vector
$\Pi_{i<j} e_{ij}^{k_{ij}} u$ has the $T_0$-weight
$$\l+\sum_{\underline{m}<i<j}k_{ij}(\e_i-\e_j)-\sum_{i\leq \underline{m}<j}k_{ij}\e_j\in
\L,$$ and the $T_0$-height $|\l|-\sum_{i\leq
\underline{m}<j}k_{ij}$.

Let $\mathcal S_l$ be the subset of $\mathcal S$ consisting of
elements with the minimal $T_0$-height $l\in \mathbf Z$, and let
$\mathcal N_l$ be the subspace of $\mathcal N$ spanned by $\mathcal
S_l$. Then it is easy to see that $e_{ij}\mathcal S_l=0$ for any
$i\leq \underline{m}<j$, and $e_{ij}\mathcal N_l\subseteq \mathcal
N_l$ for any $\underline{m}<i<j$ or $i<j\leq \underline{m}$.
\par Each element in $\mathcal S_l$ has the weight
 $$\l'+ \sum_{\underline{m}<i<j}k_{ij}(\e_i-\e_j),
\quad |\l'|=l.$$ Let $\mathcal S^h_l$ be the subset of $\mathcal
S_l$ consisting of elements
 with the largest number $\sum_{\underline{m}<i<j}k_{ij}$, and let
 $\mathcal N_l^h$ be the subspace spanned by $\mathcal S^h_l$.
  Then  we have $e_{ij}\mathcal S^h_l=0$ for all $\underline{m}<i<j$,
 and $e_{ij}\mathcal N^h_l\subseteq \mathcal N^h_l$ for all $i<j\leq \underline{m}$.
This implies that $\mathcal N^h_l$ is $u(N^+_n)$-submodule, where
$N^+_{n}$ is the sum of all positive root spaces of $\g_{n}.$
  By \cite[Th.3.2]{f}, there is a homogeneous
weight vector $0\neq v\in \mathcal N^h_l$, such that $N^+_{n}\cdot
v=0$, and hence $N^+\cdot v=0$. Thus, $v\in \mathscr M$ is a maximal
vector.\end{proof}
\par By the lemma, each simple $u_{\chi}(P_0)$-module $\mathcal M$ in
Category $\mathcal{O}$  is a homomorphic image of the induced
$u_{\chi}(P_0)$-module defined above, and hence each simple
$u_{\chi}(\g)$-module in  Category $\mathcal{O}$ is a homomorphic
image of the induced module $Z^{\chi}(\mathcal M(\l))$.\par
\subsection{The simplicity of $Z^{\chi}(\mathcal M(\l))$}
   With the same assumption as
above. Recall $\mathscr N^+_m=\oplus ^m_{i=1}\g_{-\d_i}$. Suppose
$\d_{i_1},\cdots,\d_{i_m}$ is the roots of $\mathscr N^+_m$ in the
order of increasing heights. By Th. 2.13, $u(\mathscr N^+_m)$ has a
basis
$$\{f_{\d_1}^{a_1}\cdots f_{\d_m}^{a_m}|0\leq a_i\leq \bar p-1.\}$$
For the convenience, we denote  $f_{\d_1}^{a_1}\cdots
f_{\d_m}^{a_m}\in u(\mathscr N^+_m)$ by $[a_1,a_2,\cdots,a_m]$, and
$$f^{a_1}_{\d_1}\cdots f^{a_k}_{\d_k}
f_{\d_j}f^{a_{k+1}}_{\d_{k+1}}\cdots f^{a_m}_{\d_m}\in u(\mathscr
N^+_m)$$ by $[a_1,\cdots,a_{k-1},f_{\d_j}\cdots ]$.
\begin{lemma}  Let $[a_1,\cdots, a_m]\in u(\mathscr N^+_m)$
and let $1\leq s\leq m$.  If $a_{i_j}=\bar p-1$ for all $j\geq s$,
then $[a_1,\cdots,a_{k-1},f_{\d_{i_s}},\cdots]=0$ for any $1\leq
k\leq m+1$.
\end{lemma}\begin{proof}  Let us note that if $\g_{\d}\subseteq \g_{\a}$
 and if $\a\in\G^-$, then it is possible
that $2\d$ is also a root of $\g$ (see \cite[p.51]{k3}). In this
case $2\d$ is a root of $\g^+$ with the height greater than that of
$\d$. Let $f_{\d}\in \g_{-\d}$. Then we have
$f^2_{\d}=\frac{1}{2}[f_{\d},f_{\d}]\in \g_{-2\d}$.\par In
$u(\mathscr N^+_m)$, we have
$$(*)\quad [f_{\d_i},f_{\d_j}]=f_{\d_i}f_{\d_j}-(\bar f_{\d_i},\bar f_{\d_j})f_{\d_j}f_{\d_i}=
\begin{cases}c f_{\d_i+\d_j},&\text{if $\d_i+\d_j$ is a root}\\
0,&\text{otherwise,}\end{cases}$$ where $ c\in\mathbf F^{\times}$.
We shall proceed with  induction on $s$. Suppose $s=m$. Then since
$\d_{i_m}$ has the greatest height, $[f_{\d_{i_m}}, f_{\d_j}]=0$ for
all $1\leq j\leq m$. The fact that $f^{\bar p}_{\d_{i_m}}=0$ gives
$$[a_1,\cdots,a_{k-1},f_{\d_{i_m}},\cdots]=0,$$ for any $1\leq k\leq m+1$ whenever
$a_{i_m}=\bar p-1$.\par Suppose $[a_1,\cdots,a_m]\in u(\mathscr
N^+_m)$ with $a_{i_j}=\bar p-1$ for all $j\geq s$ and $s<m$. For any
$1\leq k\leq m+1$, if $\d_s=\d_k$, then
$$[a_1,\cdots,a_{k-1},f_{\d_s},\cdots]=[a_1,\cdots,a_{k-1},\bar
p,\cdots]$$$$=\begin{cases}c[a_1,\cdots,a_{k-1},f_{2\d_k},\bar p-2,\cdots], &\text{if $2\d_k$ is a root}\\
0,&\text{otherwise.}\end{cases}$$ Since the height of $2\d_k$ is
greater than that of $\d_k$, the induction hypothesis yields
$[a_1,\cdots,a_{k-1},f_{\d_s},\cdots]=0$. The case $\d_s=\d_{k-1}$
is similar.\par  Now suppose $\d_s\notin\{\d_{k-1},\d_k\}$. By
repeated applications of the formula $(*)$ and the induction
assumption, we have
$$[a_1,\cdots,a_{k-1},f_{\d_s},\cdots]=c[a_1,\cdots,a_{s-1},f_{\d_s}
,\cdots]$$ for some $c\in\mathbf F$. Then the discussion above
applied, we get $[a_1,\cdots,a_{k-1},f_{\d_s},\cdots]=0$.

\end{proof}
It follows from the lemma that  $$f_{\d_i}[\bar p-1,\cdots,\bar
p-1]=0\quad [\bar p-1,\cdots,\bar p-1]f_{\d_i}=0,$$ for all $1\leq
i\leq m$. \par Recall the induced module $ Z^{\chi}(\mathcal
M(\l))$. Let $v\in\mathcal M(\l)$ be a maximal vector of $H$-weight
$\l\in H^*$. By applying the Harish-Chandra homomorphism, we get
$$e^{\bar{p}-1}_{\d_1}\dots
e^{\bar{p}-1}_{\d_m}f^{\bar{p}-1}_{\d_1}\cdots f^{\bar{p}-1}_{\d_m}
\otimes v=f(h)\otimes v=1\otimes f(\l(h))v,$$ where
$$f(h)=\gamma(e^{\bar{p}-1}_{\d_1}\dots
e^{\bar{p}-1}_{\d_m}f^{\bar{p}-1}_{\d_1}\cdots
f^{\bar{p}-1}_{\d_m})\in u_{\chi}(H).$$
\begin{theorem}  $Z^{\chi}(\mathcal M(\l))$ is simple if and only if
$f(\l(h))\neq 0$.
\end{theorem}
\begin{proof} Let $\mathcal N$ be a simple submodule of $Z^{\chi}(\mathcal M(\l))$.
Let us take $$0\neq v=\sum f^{k_1}_{\d_1}\cdots
f^{k_m}_{\d_m}\otimes m_{k_1,\dots,k_m}\in \mathcal N,$$ where each
$m_{k_1,\dots,k_m}\in \mathcal M(\l)$ is homogeneous. Let
$\d_{i_1},\dots,\d_{i_m}$ be in the order of increasing heights.
Using Lemma 3.10 and  applying $f_{\d_{i_m}}, \dots, f_{\d_{i_1}}$
repeatedly,    we get $0\neq f^{\bar p-1}_{\d_1}\cdots f^{\bar
p-1}_{\d_m}\otimes m\in \mathcal N$, with $m$ homogeneous. It then
follows that $$ e^{\bar p-1}_{\d_1}\cdots e^{\bar p-1}_{\d_m}f^{\bar
p-1}_{\d_1}\cdots f^{\bar p-1}_{\d_m}\otimes m\in \mathcal N.$$ By
the assumption $\chi(\g_{\pm\d_i})=0$, $1\leq i\leq m$, we can
regard $e^{\bar p-1}_{\d_1}\cdots e^{\bar p-1}_{\d_m}f^{\bar
p-1}_{\d_1}\cdots f^{\bar p-1}_{\d_m}\in u_{\chi}(\g)$ as an element
in some $u_{\chi'}(\g)$ with $\chi'$ being standard semisimple, so
that we can apply the Harish-Chandra homomorphism $\gamma$. Using
the fact that $e_{\d_i}\cdot \mathcal M(\l)=0$, $i=1,\dots, m$ and
applying $\gamma$, we get $1\otimes f(\l(h))m\in \mathcal N$. Then
the assumption $f(\l(h))\neq 0$ says that $1\otimes m\in \mathcal
N$, hence $\mathcal M(\l)\subseteq \mathcal N$. Thus, we get
$\mathcal N=Z^{\chi}(\mathcal M(\l))$, so that $Z^{\chi}(\mathcal
M(\l))$ is simple.\par
 Let $(\d_i)_{i\in I}$ be the positive root
 system of $\mathcal Z$. Then  $(\d_i)_{i\in I}\cup\{\d_1,\dots,\d_m\}$
 is a positive root system of $\g$.    For any $j\in I$, if
 $[f_{\d_{j}},f_{\d_i}]\neq 0$, then it must be a nonzero multiple
 of $f_{\d_i+\d_{j}}$.  The height of  ${\d_i+\d_{j}}\in\{\d_1,\dots,\d_m\}$ is
  greater than that of $\d_i$. Using the formula for $ad x|_{U(\g)}$ in
 2.3 and Lemma 3.10, we get  $[f_{\d_{j}}, f^{\bar
 p-1}_{\d_1}\cdots f^{\bar p-1}_{\d_m}]=0$. If
 $[e_{\d_{j}},f_{\d_i}]\neq 0$, then it must be a multiple of
 $f_{\d_i-\d_{j}}$ with  $\d_i-\d_{j }\in\{\d_1,\dots,\d_m\}$, since $\mathscr N^+_m$ is an ideal of
 $P_0^+=\mathcal Z\oplus \mathscr N^+_m$.
 Applying Lemma 3.10, we get $$[e_{\d_{j}},f^{\bar
 p-1}_{\d_1}\cdots f^{\bar p-1}_{\d_m}]=0.$$
This implies that $\mathcal M_0=f^{\bar
 p-1}_{\d_1}\cdots f^{\bar p-1}_{\d_m}\otimes \mathcal M(\l)$ is a
 $\mathcal Z$-submodule of $Z^{\chi}(\mathcal M(\l))$. By Lemma 3.10, $\mathcal M_0$ is annihilated
 by $\mathscr N^+_m$, so it is a $u_{\chi}(P_0^-)$-submodule.
 Recall  $\mathscr N^+_0=\oplus^m_{i=1}\g_{\d_i}$.
 The simplicity of $Z^{\chi}(\mathcal M(\l))$ says that
  $$u(\mathscr N_0^+)\mathcal M_0=u(\mathscr N_0^+)u_{\chi}(P_0^-)\mathcal M_0=u_{\chi}(\g)\mathcal M_0=Z^{\chi}(\mathcal M(\l)).$$  Let $(v_i)_{i\in I}$ be a basis
  of $\mathcal M_0$.  Then $$\{e^{k_1}_{\d_1}\cdots
  e^{k_m}_{\d_m}v_i|0\leq k_i\leq\bar p-1,i\in I\}$$ is a
  basis of $Z^{\chi}(\mathcal M(\l))$. Taking  $f^{\bar
 p-1}_{\d_1}\cdots f^{\bar p-1}_{\d_m} \otimes v\in\mathcal M_0$, where $0\neq v\in \mathcal M(\l)$
 is a maximal vector of $H$-weight $\l$,
  we get $$0\neq e^{\bar p-1}_{\d_1}\cdots e^{\bar p-1}_{\d_m}f^{\bar
 p-1}_{\d_1}\cdots f^{\bar p-1}_{\d_m}\otimes v=1\otimes f(\l(h))v,$$ so that
 $f(\l(h))\neq 0$.
\end{proof}
\subsection{The determination of $f(\l(h))$}
Recall $\mathscr N_k^+=\sum^{k}_{i=1}\g_{-\d_i}+\sum^m_{i=k+1}
\g_{\d_i}$, $0\leq k\leq m$.
\begin{lemma}
With the same assumption as in 3.6. For each $0\leq k\leq m$, we
define a subspace of $Z^{\chi}(\mathcal M(\l))$ by
$$M_k=:f^{\bar p-1}_{\d_k}\cdots f^{\bar p-1}_{\d_1}\otimes
\mathcal M(\l). $$  If $2\d\notin \{\d_1,\dots,\d_m\}$ for any
$\d\in\{\d_1,\dots,\d_m\}$, then $M_k$ is annihilated by $\mathscr
N^+_k$.

\end{lemma}\begin{proof}We proceed with induction on $k$. The case $k=0$ is
given by the definition. Assume the claim is true for $M_{k-1}$ with
$k\geq 1$.
 By definition, $M_k=f^{\bar p-1}_{\d_k}M_{k-1}$. Using Lemma 3.10 and the
 assumption  that $\{-\d_1,\dots,
 -\d_k\}$ is a closed system, we get $f_{\d_i}M_k=0$ for all
 $i=1,\dots,k$.\par
  Suppose $\d=\d_i$ and $i\geq k+1$.  Then $\d$ is also a positive root of $\mathscr N^+_{k-1}$
 and $\d\neq\d_k$, so it suffices to assume $[e_{\d},f_{\d_k}]\neq 0$. Since $\d_k$ is simple in the
 positive root system $\Phi^+_{k-1}$,  $\d-\d_k$ a positive root
 for $\Phi^+_{k-1}$.  Then we must have
  $[e_{\d},f_{\d_k}]=ce_{\d-\d_k}$ for some
 $c\in\mathbf F^{\times}$.  \par If $\d_k$ is root for $\g^+$, by
 assumption we get
 $\d-\d_k\in\Phi^+_{k-1}\setminus \{\d_k\}$. Then since $\d_k$ is simple,
 we have for any $i\in\mathbf N$, $\d-i\d_k\in\Phi^+_{k+1}\setminus \{\d_k\}$ if
 $\d-i\d_k\in\Phi$,  and so the induction
assumption yields
 $$e_{\d}f^{\bar p-1}_{\d_k}M_{k-1}=([e_{\d},f_{\d_k}]f^{p-2}_{\d_k}+
 \dots +f^{p-2}_{\d_k}[e_{\d},f_{\d_k}])M_{k-1}$$$$=c(e_{\d-\d_k}f_{\d_k}^{p-2}+\dots +
 f^{p-2}_{\d_k}e_{\d-\d_k})M_{k-1}=0.$$\par If $\d_k$ is a root for
 $\g^-$, then we have $\bar p=1$, so that $$e_{\d}\cdot
 M_k=e_{\d}f_{\d_k}M_{k-1}$$$$=(\bar e_{\d},\bar f_{\d_k})f_{\d_k}e_{\d}
 M_{k-1}+ce_{\d-\d_k}M_{k-1}=0,$$ where the last equality  is given by induction
 hypothesis.
\end{proof}

 Let $A=\oplus _{\a\in\G}A_{\a}$ be
a finite dimensional Frobenious algebra. We use the notation
$A_L$(resp. $A_R$) to denote the left(resp. right) $\G$-graded
regular $A$-module $A$. For each $\a\in \G$, let
$$A_{\a}^*=\{f\in A^*|f(A_{\b})=0, \text{for all}\quad \b\neq
\a\}.$$ Then we have $A^*=\oplus _{\a\in\G} A^*_{\a}$. It is easy to
check that $A_L^*$(resp. $A_R^*$) is a right(resp. left) $\G$-graded
$A$-module.
\begin{lemma}\cite[Th. 61.3]{cr} For a finite dimensional
$\G$-graded algebra $A$, the following are equivalent.\par (1) $A$
is Frobenious.\par (2) $A_L\cong A^*_R$, $A_R\cong A_L^*$.

\end{lemma}
\begin{definition}\cite{f}Let $\g$ be a restricted Lie color algebra. We call
$\g$ unipotent if for every $x\in \g_{\a}$ with $\a\in \G^+$, there
exists $r>0$ such that $x^{[p]^r}=0$.
\end{definition}

\begin{proposition}Let $\g=\oplus_{\a\in\G}\g_{\a}$ be a finite dimensional unipotent Lie
color algebra. Then the left(resp. right) regular $u(\g)$-module
$u(\g)$ has a unique 1-dimensional trivial submodule.
\end{proposition}
\begin{proof}   Let $x\in\g_{\a}$ with
 $\a\in\G^+$. Then we have $(x^p-x^{[p]})u(\g)=0$.  Since $\g$ is
unipotent, there is $n\in \mathbf Z^+$ such that $x^{[p]^n}=0$.
Hence we get $x^{p^n}u(\g)=0$. If $x\in\g_{\a}$, $\a\in\G^-$, then
we get
 $x^2=\frac{1}{2}[x,x]
\in \g_{2\a}$,  so  $x$ is also nilpotent in $u(\g)$, since
$2\a\in\G^+$. By \cite[Th. 3.2]{f}, $u(\g)$ has a  1-dimensional
trivial submodule $\mathbf Fv$. \par To show the uniqueness, we let
$f$ be the nondegenerate invariant bilinear form on $u(\g)$. Then we
get $f(x,v)=f(1,xv)=0$ for all homogeneous $x\in \g$, so that $v$ is
in the homogeneous right orthogonal complement of $u(\g)\g$. Since
$u(\g)\g$ has codimension 1 in $u(\g)\g$, its homogeneous orthogonal
complement is 1-dimensional. This implies that $u(\g)$ has a unique
trivial submodule $\mathbf Fv$.
\end{proof}

\begin{proposition}\cite{j} Let $\g$ be a finite dimensional unipotent Lie color
algebra, and let $\chi=\{\chi_{\a}|\a\in \G^+\}$ be a $p$-character
of $\g$.  Then  $u_{\chi}(\g)$ has only one(up to isomorphism)
simple module.
\end{proposition}
\begin{proof}Let $M=\oplus _{\a\in\G}M_{\a}$ and
$M'=\oplus_{\a\in\G}M'_{\a}$ be two simple $u_{\chi}(\g)$-modules.
By \cite[Lemma 2.5]{bm}, we have $\dim M,\dim M'<\infty$. For any
$\a\in\G$, let
$$\text{Hom}_{\mathbf F}(M,M')_{\a}=\{f\in \text{Hom}_{\mathbf F}(M,M')|f(M_{\b})\subseteq
M'_{\b+\a}\}.$$ Then $\text{Hom}_{\mathbf F}(M,M')=\oplus_{\a\in\G}
\text{Hom}_{\mathbf F}(M,M')_{\a}$. $\text{Hom}_{\mathbf F}(M,M')$
is a  $\g$-module with the $\g$-action defined by(\cite{f})
$$(x\cdot f)(m)=x\cdot f(m)-(\bar x|\bar f)f(x\cdot m)$$ for all
homogeneous $x$ and $f$. Then one can easily check that
$\text{Hom}_{\mathbf F}(M,M')$ is a restricted $\g$-module.\par For
each $x\in\g_{\a}$, by the proof of Proposition 3.15, we see that
$x$ acts nilpotently on $\text{Hom}_{\mathbf F}(M,M')$. A
straightforward computation shows that $\text{Hom}_{\mathbf
F}(M,M')_{\bar 0}$ is a $u(\g)$-submodule of $\text{Hom}_{\mathbf
F}(M,M')$.  By \cite[Th.3.2]{f}, $\text{Hom}_{\mathbf F}(M,M')_{\bar
0}$ contains an 1-dimensional trivial submodule $\mathbf Ff$.  By
definition,   $f$ is a homomorphism of $\g$-modules. Then the
simplicity of both $M$ and $M'$ implies that $f$ is an isomorphism.
Thus,  $M\cong M'$.
\end{proof}

\begin{lemma}Let $\mathbf Fv_L$(resp. $\mathbf Fv_R$) be the unique trivial
$u(\g)$-submodule of $u(\g)_L$(resp. $u(\g)_R$). Then $ v_L=c v_R$,
for some $c\in\mathbf F^{\times}$.
\end{lemma}
\begin{proof}By Lemma 2.19, $u(\g)$ is a symmetric algebra. Let $f$ be
a nondegenerate invariant bilinear form on $u(\g)$. The symmetry of
$f$ implies that both $v_L$ and $v_R$ are in the right homogeneous
orthogonal complement of $u(\g)\g$, which is 1-dimensional. Then we
get $ v_L=c v_R$ for some $c\in\mathbf F^{\times}$.
\end{proof}
Let $(\g, Z,\{\d_1,\dots,\d_m\})$ be a FP triple. Recall $\mathscr
N^+_m=\oplus^m_{i=1}\g_{-\d_i}$.  By the remark following Lemma
3.10,  $\mathbf Ff_{\d_1}^{\bar p-1}\cdots f_{\d_n}^{\bar
p-1}$(resp.  $\mathbf F f_{\d_n}^{\bar p-1}\cdots f_{\d_1}^{\bar
p-1}$) is a 1-dimensional trivial submodule of the $\G$-graded
regular $u(\mathscr N^+_m)$-module $u(\mathscr N^+_m)_L$(resp.
$u(\mathscr N^+_m)_R$), so we get
$$f_{\d_1}^{\bar p-1}\cdots f_{\d_n}^{\bar p-1}=cf_{\d_n}^{\bar
p-1}\cdots f_{\d_1}^{\bar p-1}$$ for some $c\in \mathbf F^{\times}$.
\begin{lemma}Let $\g=\oplus_{\a\in\G}\g_{\a}$ be a Lie color
algebra, and let $M=\oplus_{\a\in\G}M_{\a}$ be a $\g$-module. We
take $e\in\g_{\a}$,  $f\in\g_{-\a}$ and assume $[e,f]=h$. Let $v\in
M$ such that $e\cdot v=0$.
\par (1) If $(\a,\a)=-1$, then $ef v=(fe+h)v$.
\par (2) If $(\a|\a)=1$, then $$e^lf^l
v=l!\Pi^{l-1}_{i=0}(h-i) v$$ for every $l>1$.\par (3) If $h\cdot
v=\l(h) v$, for some $\l(h)\in\mathbf F$, then $$e^{\bar p-1}f^{\bar
p-1}v=(\bar p-1)![(\l(h)+1)^{\bar p-1}-1]v.$$
\end{lemma}
\begin{proof}(1) is obvious. (2) is given by induction on $l$. \par (3)
From (1) and (2), we get$$ e^{\bar p-1}f^{\bar p-1}v=(\bar
p-1)!\Pi^{\bar p-2}_{i=0}(\l(h)-i)v.$$ Then the claim follows from
the identity $$ x^{p-1}-1=\Pi^{p-1}_{i=1}(x-i).$$
\end{proof}
\begin{theorem} Let $f(\l(h))$ be the polynomial
  given in Sec. 3.6. Assume $2\d\notin\{\d_1,\dots,\d_m\}$, for any $\d\in\{
  \d_1,\dots,\d_m\}$. Then  we have
$$f(\l(h))=c\Pi^m_{i=1}[(\l_i(H_{\d_i})+1)^{\bar p-1}-1]$$ for some
$c\in \mathbf F^{\times}$,  where $\l_i=\l-[(\bar p-1)\d_1+\dots
+(\bar p-1)\d_{i-1}]$, $1\leq i\leq m$. \end{theorem} \begin{proof}
Since $f^{\bar p-1}_{\d_m}\cdots f^{\bar p-1}_{\d_1}=kf_{\d_1}^{\bar
p-1}\cdots f_{\d_m}^{\bar p-1}$ in $u(\mathscr N^+_m)$,  $k\in
\mathbf F^{\times}$, we get
$$e^{\bar p-1}_{\d_1}\cdots e_{\d_m}^{\bar p-1}f_{\d_1}^{\bar
p-1}\cdots f_{\d_m}^{\bar p-1}\otimes v$$$$=k e^{\bar
p-1}_{\d_1}\cdots e_{\d_m}^{\bar p-1}f_{\d_m}^{\bar p-1}\cdots
f_{\d_1}^{\bar p-1}\otimes v$$$$=ke_{\d_1}^{\bar p-1}\cdots
e_{\d_{m-1}}^{\bar p-1}(e^{\bar p-1}_{\d_m}f_{\d_m}^{\bar
p-1})f^{\bar p-1}_{\d_{m-1}}\cdots f_{\d_1}^{\bar p-1}\otimes
v$$(using Lemma 3.12 and Lemma 3.18)$$=k(\bar
p-1)![(\l_m(H_{\d_m})+1)^{\bar p-1}-1]e_{\d_1}^{\bar p-1}\cdots
e_{\d_{m-1}}^{\bar p-1}f^{\bar p-1}_{\d_{m-1}}\cdots f_{\d_1}^{\bar
p-1}\otimes v$$$$=\dots =c\Pi^m_{i=1}[\l_i(H_{\d_i})+1)^{\bar
p-1}-1]v,$$  $ c\in \mathbf F^{\times}$.  Thus, the claim holds.
\end{proof}

\section{Applications}
In this section, we consider the applications of Th. 3.11 and Th.
3.19 to the Lie color algebra $\g=\text{cgl}(V)$.   We prove an
analogue of Kac-Weisfeiler theorem.  Then we determine the condition
for the baby Verma module to be simple.
\subsection{The Kac-Weisferlar theorem} Let $\g=\text{cgl}(V)$. The $\mathbf F$-vector
space $\text{cgl}(V)$ also has a natural $\G$-graded restricted Lie
algebra structure with the Lie product defined by $[x,y]=xy-yx$,
$x,y\in \cup_{\a\in\G}\text{cgl}(V)_{\a}$. We denote this Lie
algebra by $\text{cgl}(V)^-$. Recall that $\g$ has a filtration
$$\g_1\subseteq \g_2\subseteq \dots \subseteq
\cup_{i=1}^{\infty}\g_i=\g$$ in case $\g$ is infinite dimensional.
Then the underlining $\mathbf F$-vector subspace of each $\g_n$
becomes a Lie subalgebra of $\text{cgl}(V)^-$,  denoted
$\g_n^-$.\par

 Recall the linear algebraic group $G$
with Lie$(G)=\g_{\bar 0}$. Then $G$ is also a group of automorphisms
for $\text{cgl}(V)^-$ which keeps both the $\G$-grading and the
$p$-mapping. Moreover, $\text{cgl}(V)^-$ has the same root space
decomposition relative to $T$ with that of $\text{cgl}(V)$. i.e.,
$\Phi^+$(resp. $\Delta$) in 3.1 is the positive(resp. simple) root
system for both $\text{cgl}(V)^-$ and $\text{cgl}(V)$. Recall the
notion $e_{\d}$, $f_{\d}$ for $\d\in\Phi^+$.
\par We use $[, ]_{\G}$(resp. $[,]$) momentarily to denote the Lie
color product(resp. Lie product) in $\text{cgl}(V)$. Then it is easy
to see that
$$[e_{\d_i},e_{\d_j}]_{\G}=c_{ij}[e_{\d_i},e_{\d_j}], [f_{\d_i},f_{\d_j}]_{\G}
=c_{ij}[f_{\d_i},f_{\d_j}]\quad\text{for any}\quad
\d_i,\d_j\in\Phi^+$$ and $$
[e_{\d_i},f_{\d_j}]_{\G}=k_{ij}[e_{\d_i},f_{\d_j}]
\quad\text{whenever}\quad \d_i\neq \d_j,$$ where  $c_{ij}$,
$k_{ij}\in\mathbf F^{\times}$. \subsubsection{The refections} Recall
 the simple roots $\Delta$ for $\Phi^+$:
$$ \Delta=\{\e_1-\e_2,\e_2-\e_3,\dots,\e_{n-1}-\e_n,\cdots\}.$$
 We define a real vector space $V$
having a basis $\Delta$. Then $V$ is infinite dimensional if $\G$ is
infinite. For any $\d_i,\d_j\in\Delta$, we let $$\langle
\d_i,\d_j\rangle=\begin{cases}2, &\text{if}\quad
i=j\\-1,&\text{if}\quad
i\in\{j-1,j+1\}\\0,&\text{otherwise}.\end{cases}$$ For each $\d\in
\Delta$,  we define a linear function $\langle -,\d\rangle\in V^*$
by
$$\langle \beta,\d\rangle=\sum_{i}c_i\langle \d_i,\d\rangle
\quad\text{if}\quad \beta=\sum_{i}c_i\d_i\in V.$$ For each $\d\in
\Delta$, we define a reflection on $V$ by $$\tau_{\d}
(\beta)=\beta-\langle \beta,\d\rangle \d.$$   It is easy to see that
$\tau_{\d}(\d)=-\d$, and $\tau_{\d}(\beta)=0$ for any $\beta\in V$
with $\langle \beta,\d\rangle=0$. The group $W$ generated by all
$\tau_{\d}$'s,$\d\in\Delta^+$ is called the Weyl group  of the Lie
algebra $\text{cgl}(V)^-$.  Then $W$ is an infinite group in case
$\G$ is infinite.\par Let $\g=\text{cgl}(V)$ be infinite
dimensional. Then the Lie algebra $\text{cgl}(V)^-$ has a filtration
$$\g^-_1\subseteq \g^-_2\subseteq \dots \subseteq
\cup^{\infty}_{i=1}\g^-_i=\text{cgl}(V)^-.$$ Let us denote the root
system(resp. positive root system, simple root system, Weyl group)
of $\g_n$ by $\Phi_n$(resp. $\Phi^+_n$, $\Delta_n$, $W_n$). Let
$V_n$(resp. V) denote the real vector space spanned by
$\Phi^+_n$(resp. $\Phi$). Then we get
$$\Phi=\cup^{\infty}_{n=1}\Phi_n,\quad
\Phi^+=\cup^{\infty}_{n=1}\Phi^+_n\quad
\Delta=\cup^{\infty}_{n=1}\Delta_n.$$ By extending the action of
each $\sigma\in W_n$ on $V_n$ to that of $V$, we can identify $W_n$
as a subgroup of $W$. Then we get  $W=\cup^{\infty}_{n=1}W_n$.
\par Let $\sigma\in W$ and  assume $\sigma=\tau_{\d_1}\cdots
\tau_{\d_l}$, $\d_1,\dots,\d_l\in \Delta$. Let $N$ be large enough
that $\Delta_N$ contains each $\d_i$. Then for each $n\geq N$, we
have $\sigma_{|V_n}\in W_n$, and hence $\sigma (\Phi_n)\subseteq
\Phi_n$. This gives $\sigma (\Phi)=\Phi$, and hence $\sigma(\Phi^+)$
is a system of positive roots of $\text{cgl}(V)^-$ with the set of
simple roots $\sigma(\Delta)$.

\begin{lemma}For each $\sigma\in W$, $\sigma(\Phi^+)$ is a system of
positive roots of the Lie color algebra $\g=\text{cgl}(V)$ with the
set of simple roots $\sigma(\Delta)$.\end{lemma}\begin{proof}Since
$\sigma(\Phi^+)$ is a positive root system of the Lie algebra
$\text{cgl}(V)^-$ with simple roots $\sigma(\Delta)$, $\sigma
(\Delta)$ is a minimal subset $\sigma(\Phi^+)$ satisfying:\par (1)
$\{e_{\gamma}|\gamma\in \sigma(\Delta)\}$ (resp.
$\{f_{\gamma}|\gamma\in s_{\a}(\Delta)\})$ generates
$\{e_{\gamma}|\gamma\in \sigma(\Phi^+)\}$(resp.
$\{f_{\gamma}|\gamma\in \sigma(\Phi^+)\}$).\par
(2)$\{e_{\gamma},f_{\gamma}|\gamma\in \sigma(\Delta)\}$ generates
$\text{cgl}(V)^-$.\par (3) $[e_{\a},f_{\b}]=\delta_{\a,\b}h\in H$,
for any $\a$,$\b\in \sigma(\Delta)$.\par Since each $e_{\gamma}$
above is also a root vector for the Lie color algebra
$\g=\text{cgl}(V)$, the formulas preceding 4.1.1 implies that
$\sigma(\Phi^+)$ is a positive root system of $\text{cgl}(V)$ with
the set of simple roots $\sigma(\Delta)$.
\end{proof}Let $\g=\text{cgl}(V)$ and let $\chi\in\g^*_{\bar 0,f}$
be in the standard form. We consider both the case $\g$ is finite
dimensional and the case $\g$ is infinite dimensional.   Let
$$\mathcal Z=c_{\g}(\chi_s)=\oplus_{\a\in\G}\mathcal Z_{\a},$$ where
$\mathcal Z_{\a}=\{x\in \g_{\a}|\chi_s([x,-])=0\}.$ Then by Example
(1), (3) of 3.4, we get a resulted FP triple $(\g, \mathcal Z,
(\d_i)^m_{i=1})$.  Recall the parabolic color subalgebra ${P_0}$
having
$$\mathscr N_0^+=\sum_{\d\in\Phi^+,\chi (H_{\a})\neq 0}\g_{\d}=\oplus_{i=1}^m\g_{\d_i}$$ as its nilradical.
\par
Recall the definition of $Z^{\chi}(\mathcal M(\l))$ and $f(\l(h))$.
Since $2\d\notin \Phi^+$ for every $\d\in\Phi^+$, we get by Th. 3.19
that
$$f(\l(h))=c\Pi^m_{i=1}[(\l_i(H_{\d_i})+1)^{\bar p-1}-1],\quad c\in\mathbf F^{\times}.$$
For each $1\leq i\leq m$, we have $H^{[p]}_{\d_i}=H_{\d_i}$.  This
implies that $$\l(H_{\d_i})^p-\l(H_{\d_i})=\chi(H_{\d_i})^p\neq 0,$$
and hence $\l(H_{\d_i})\notin \mathbf F_p$. Thus  $f(\l(h))\neq 0$.
Then we get from Th.3.11 that
\begin{corollary}(Kac-Weisfeiler Theorem) Let $\g=\text{cgl}(V)$. With the assumption
as above, then $Z^{\chi}(\mathcal M(\l))$ is a simple
$u_{\chi}(\g)$-module.
\end{corollary}
\subsection{The simplicity of the baby
Verma module}In this subsection, we consider the simplicity of the
baby Verma module for the algebraic Lie color algebra
$\g=\text{cgl}(V)$. \subsubsection{$\g$ is finite dimensional} Let
$\G$ be finite. Then $\text{cgl}(V)$ is the general linear Lie color
algebra $\g=\text{gl}(m,\G)$.  Recall the  triangular decomposition
of $\g$: $\g=N^+\oplus H\oplus N^-$.  We assume $\chi\in\g^*_{\bar
0}$ is standard semisimple.  Then by Example (2) of 3.4, there is a
resulted FP triple: $(\g,H,\Phi^+)$.\par In the present situation,
the parabolic subalgebra $P_0$ is just the Borel subalgebra
$B=H\oplus N^+$, and $\mathscr N^+_0=N^+$ is its nilradical. By
\cite[Th. 3.2]{f}, each simple $u_{\chi}(B)$-module $\mathcal M(\l)$
is 1-dimensional and defined as follows:\par
$$\mathcal M(\l)=\mathbf F v, N^+ v=0, h\cdot v=\l(h) v \text{ for every} \quad h\in
H,$$ where $\l(h)\in H^*$ satisfies
$\l(H_{\d_i})^p-\l(H_{\d_i})=\chi^p(H_{\d_i})$ for any
$\d_i\in\Phi^+$.  Recall the induced module $Z^{\chi}(\mathcal
M(\l))=: u_{\chi}(\g)\otimes _{u_{\chi}(B^+)} \mathbf Fv$. We call
it the baby Verma module with character $\chi$ and denote it by
$Z^{\chi}(\l)$. Then by Th. 3.11, we get
\begin{corollary} $Z^{\chi}(\l)$ is simple if and only if
$f(\l)\neq 0$.
\end{corollary}Note: In the special case that $\g$ is a finite dimensional semisimple Lie algebra,
a conclusion analogous to  Corollary 4.3 was proved by Rudakov in a
 different approach \cite{r}.\par By Th. 3.19, we get
$f(\l)=c\Pi_{\d_i\in\Phi^+}[(\l_i(H_{\d_i})+1)^{\bar p-1}-1],
c\in\mathbf F^{\times}.$ This formula enables us to conclude with
\begin{corollary}Let $\g=\text{gl}(m,\G)$ and assume that $\chi\in\g^*_{\bar 0}$ is regular semisimple,
that is,  $\chi(H_{\d})\neq 0$ for every $\d\in\Phi^+$. Then the
baby Verma module $Z^{\chi}(\l)$ is simple.
\end{corollary}

\subsubsection{$\g$ is infinite dimensional}
We now let $\g=\text{cgl}(V)$ be infinite dimensional. Recall the
triangular decomposition of $\g$ in 3.1: $\g=N^+\oplus H\oplus N^-$.
Let $B$ denote the Borel subalgebra $H\oplus N^+$.  Assume
$\chi\in\g^*_{\bar 0,f}$ is standard semisimple. We take $T_0$ as
that in the Example 3.2. By Lemma 3.9, any simple $u_{\chi}(B)-T_0$
module in Category $\mathcal{O}$ is $1$-dimensional. Let us denote
it by $\mathbf Fv$ and assume $\l\in H^*$ is the $H$-weight of $v$.
Then the ($\G$-graded)baby Verma module is defined as
$$Z^{\chi}(\l)=u_{\chi}(\g)\otimes _{u_{\chi}(B)}\mathbf Fv.$$ It is
easy to see that $Z^{\chi}(\l)$, as a $u_{\chi}(\g)-T_0$ module, is
also an object in Category $\mathcal{O}$. \par  Recall the finite
dimensional Lie color subalgebra $\g_n$ for each $n\in\mathbf Z^+$.
Let $H_n$ be its maximal torus consisting of diagonal matrices, and
let $\d_1,\dots,\d_{k_n}$ be the set of all its positive roots. By
the arguments in \cite[Lemma 8.4]{f1} and Lemma 4.1, the set of
roots can be ordered such that ($\g_n,H_n, \d_1,\dots,\d_{k_n})$ is
a FP triple. Then we obtain a polynomial:$$
 f_n(\l(h))=\Pi^{k_n}_{i=1}[\l_i(H_{\d_i})+1)^{\bar p-1}-1].$$
\begin{theorem}  $Z^{\chi}(\l)$ is simple if and only if $f_n(\l(h))\neq 0$ for
all $n\geq 1$.
\end{theorem}
\begin{proof}Assume $f_n(\l(h))\neq 0$ for
all $n\geq 1$. Let $\mathcal N\subseteq Z^{\chi}(\l)$ be a simple
$u_{\chi}(\g)$-submodule. Let us take an element $$0\neq \sum
f^{l_1}_{\d_1}\cdots f^{l_k}_{\d_k}\otimes v\in \mathcal N,\quad
0\leq l_i\leq \bar p.$$ Let
 $n$ be large enough that the root system of $\g_n$
contains all $\d_i's$ appeared in the summation. Then by applying a
similar argument as that used in the proof of Th. 3.11, we get
$1\otimes f_n(\l(h)) v\in \mathcal N$. This gives $1\otimes v\in
\mathcal N$, and hence $\mathcal N=Z^{\chi}(\l)$. Thus,
$Z^{\chi}(\l)$ is simple.\par On the other hand, suppose
$f_n(\l(h))=0$ for some $n$. Then by Coro. 4.3, the baby Verma
module $Z_n^{\chi}(\l)$ for the Lie color subalgebra $\g_n$ fails to
be simple. Here $\chi$(resp. $\l$) is the restriction of that for
$\g$(resp. $H$) to $\g_n$(resp. $H_n$). Let us take the parabolic
subalgebra $\mathcal P_n=\g_n+B$. We can write $\mathcal P_n$ in the
form
$$\mathcal P_n=\g_n\oplus H^c_n\oplus N^{(n)},$$ where $N^{(n)}$ is
the nilradical of $\mathcal P_n$ and $$H^c_n=\sum _{l>\sum_{i\leq n}
m_{\a_i}}\mathbf Fe_{ll}.$$ Since $[H^c_n,\g_n]=0$ and $N^{(n)}$
annihilates the maximal vector $v$ of $Z^{\chi}(\l)$, $N^{(n)}$
annihilates the induced $u_{\chi}(\mathcal P_n)$-submodule
$$u_{\chi}(\mathcal P_n)\otimes _{u_{\chi}(B)}\mathbf Fv\subseteq
Z^{\chi}(\l).$$  Also $H^c_n$ acts as scalar multiplications on
$u_{\chi}(\mathcal P_n)\otimes _{u_{\chi}(B)}\mathbf Fv$, so that
$u_{\chi}(\mathcal P_n)\otimes _{u_{\chi}(B)}\mathbf Fv$ is
isomorphic to $Z_n^{\chi}(\l)$ as a $u_{\chi}(\g_n)$-module.
Conversely, we can regard each $u_{\chi}(\g_n)$-module
$Z^{\chi}_n(\l)$ as the induced $u_{\chi}(\mathcal P_n)$-module
above by letting $N^{(n)}$ annihilate $Z^{\chi}_n(\l)$, and letting
$H^c_n$ act on $Z^{\chi}_n(\l)$ as  multiplications by $\l(h)$,
$h\in H^c_n$. Therefore, there is an isomorphism of
$u_{\chi}(\g)$-modules:
$$\Psi: u_{\chi}(\g)\otimes_{u_{\chi}(\mathcal
P_n)}Z^{\chi}_n(\l)\longrightarrow Z^{\chi}(\l).$$ Since
$u_{\chi}(\g)\otimes_{u_{\chi}(\mathcal P_n)}-$ is exact,
$Z^{\chi}(\l)$ is not simple. This completes the proof.
\end{proof}
\section{Appendix: The linear algebraic group
$\text{GL}(\{m_i\},\mathbf F)$ and its Lie algebra} In this section,
we define an infinite dimensional algebraic group
$\text{GL}(\{m_i\},\mathbf F)$ and its Lie algebra. We draw most of
the notation and standard procedure from \cite{hh2}.\par Let $K$ be
an algebraically closed filed. The set $K\times \cdots \times
K\cdots $ is called an infinite affine space and denoted
$A^{\omega}$. Then each element $a\in A^{\omega}$ is an infinite
sequence $a=(a_i)^{\infty}_{i=1}$, denoted  simply by $(a_i)$ in the
following. \par Let $K[T]=K[T_i]^{\infty}_{i=1}$. For each ideal $I$
in $K[T]$, let $V(I)$ denote the common zeros of all $f\in I$. Then
the collection of all $V(I)'s$ defines the Zariski topology on
$A^{\omega}$. We call each open or closed set $X$ with coordinate
ring $K[X]$ an affine variety. Similarly one can define the
morphisms of affine varieties.\par

Let $(a,b)$ denote the countable set $(a_i,b_j)^{\infty}_{i,j=1}$
for each $a=(a_i)$, $b=(b_j)\in A^{\omega}$. We define
$$A^{\omega}\times
A^{\omega}=\{(a, b)|a,b\in A^{\omega}\}.$$ The set $A^{\omega}\times
A^{\omega}$ has the Zariski topology defined with the coordinate
ring $$K[T_i,U_j]^{\infty}_{i,j=1}\cong K[T]\otimes K[U].$$ Then
$A^{\omega}\times A^{\omega}$ becomes the product of $A^{\omega}$
with $A^{\omega}$.\par Let $\{m_i|i\in \mathbf Z^+\}$ be a sequence
of positive integers, and let $\text{GL}(m_i)$ be the general linear
group consisting of invertible $m_i$ by $m_i$ matrices. We denote
the direct product
$$G=:\text{GL}(\{m_i\},\mathbf F)=\{f:\mathbf Z^+\longrightarrow
\cup^{\infty}_{i=1}\text{GL}(m_i)| f(i)\in \text{GL}(m_i)\}.$$ We
regard each element $g\in G$ as an infinite diagonal block matrix,
with the $i$th block  in $\text{GL}(m_i)$. Let $I_i$ be the set of
all pairs of integers $(i,j)$ such that $\sum^{i-1}_{k=1} m_k+1\leq
i,j\leq \sum ^i_{k=1}m_k$. Then
$$GL(m_i)=\{(a_{st})_{m_i\times m_i}|(s,t)\in I_i, \text{det} (a_{st})\neq 0\}.$$
Since the set $\cup^{\infty}_{i=1}I_i$ is countable, we can identify
$G$ with an open subset of $A^{\omega}$:
$$G=A^{\omega}-V(\{det(T_{ij})_{(i,j)\in I_i}|i\in \mathbf Z^+\}).$$\par Let
$\{A_i|i\in \mathbf Z^+\}$ be a sequence of finitely-generated
commutative $K$-algebras. We define the infinite tensor product by
$$\otimes_i A_i=:\{\otimes_i a_i|a_i=1\quad\text{for all but
finitely many i's}\}.$$ Then we can write each element
$\otimes_{i}a_i$ by a finite product, say, $a_{i_1}\otimes
\cdots\otimes a_{i_n}$ if all $a_i=1$ for $i\notin\{i_1,\dots,
i_n\}$. It follows that $\otimes_{i}A_i$ is a commutative algebra
with a countable set of generators consisting of  finite products.
Then the coordinate ring of the group $G$ is
$$K[G]=\otimes _{i}K[GL(m_i)]=\otimes_{i}
(K[T_{ij}]_{(i,j)\in I_i})_{det(T_{ij})}.$$  Let us denote the
coordinate ring of $K[GL(m_i)]$ by $K[t]_i=:K[t_{sr}|_{(s,r)\in
I_i}]$. We identify it with its canonical image in $K[G]$. By
identifying $G\times G$ with an open subset of $A^{\omega}\times
A^{\omega}$ with induced topology, we have that the group
multiplication $\phi: G\times G\longrightarrow G$ defined by
$$\phi((A_i), (B_i))=(A_iB_i)$$ is a morphism of
affine varieties, since $\phi=(\phi_{ij})_{(i,j)\in \cup_l I_l}$,
where for each $(i,j)\in I_l$, we have $$\phi_{ij}=\sum
^{m_i}_{k=1}t_{ik}u_{kj}\in K[t]_l\otimes K[u]_l.$$ Similarly, one
can show that the map $\mathfrak i:G\longrightarrow G$ defined by $
\mathfrak i(x)=x^{-1}$ is also a morphism of affine varieties. Thus,
$G$ is an algebraic group in the sense of \cite{hh2}. One can also
define as in \cite{hh2} the $G$-action on $K[G]$ via left(resp.
right) translation $\l_x$(resp. $\rho_x$).\par We are now ready to
define the Lie algebra of $G$. Let $A=K[G]$. Then $\text{Der} A$ is
a Lie algebra. We define the Lie subalgebra
$$\text{Lie}(G)=\{\d\in \text{Der} A|\d\l_x=\l_x\d, \quad \text{for
all}\quad x\in G, $$$$\d (K[\text{GL}(m_i)])=0, \text{for all i
greater than some}\quad N\in \mathbf Z^+\},$$ and call it the Lie
algebra of $G$. There are also equivalent definitions of the Lie
algebra Lie($G$) as in \cite{hh2}. For each $x\in G$, we can write
$x$ as a sequence $(x_i)$ with each $x_i\in GL(m_i)$. Let $\mathscr
O_{x_i}$ denote the local ring of the finite dimensional group
$GL(m_i)$ at $x_i$ with the unique maximal $m_{x_i}$. Let us denote
the ideal of $\otimes _i\mathscr O_{x_i}$
$$m^{(i)}_x=:\sum_{i=1}^{\infty}\mathscr O_{x_1}\otimes \cdots\otimes
\mathscr{O}_{x_{i-1}}\otimes m_{x_i}\otimes
(\otimes^{\infty}_{j=i+1}\mathscr{O}_{x_j}).$$ Then the local ring
$\mathscr O_x$ is $\otimes _i\mathscr O_{x_i}$ localized by
$m_x=\sum^{\infty}_{i=1}m_x^{(i)}$.\par Identifying each
$\mathscr{O}_{x_i}$ with its canonical image in $\mathscr{O}_x$, we
can define  $\text{Lie}(G)$  to be the set of all point derivations
from $\mathscr O_x$ to $\mathbf F$ such that $\d (\mathscr
O_{x_i})=0$ for all $i$ greater than some $N\in \mathbf Z^+$. While
in terms of tangent spaces, $\text{Lie}(G)$ can be  defined as
$$(m_x/m_x^2)^*_f=:\{\phi\in m_x^*|\phi(m^2_x)=0, $$$$\phi(m^{(i)}_x)=0, \text{for all i greater than
some}\quad N\in \mathbf Z^+ .\}$$ By  similar arguments as those
used in \cite{hh2}, one can show that all these definitions agree.
\par   Then we can define the differentiation for each
morphism  of algebraic  groups. It is easy  to see that
$\text{Lie}(G)=\oplus^{\infty}_{i=1}\text{gl}(m_i,\mathbf F)$. Each
element $x\in \text{Lie}(G)$ is a diagonal block matrix. Taking the
differentiation of the automorphism $\text{Int} x$ of $G$, $x\in G$,
one gets the adjoint action of  $G$  on $\text{Lie}(G)$:
$$\text{Ad} g(x)=gxg^{-1}, \quad x\in\text{Lie}(G), g\in G,$$ where the right
side of the equality is the product of infinite block matrices.\par
\def\refname{

\end{document}